\documentclass[]{amsart}
\usepackage{amsmath,amssymb,amsbsy,latexsym}
\usepackage{graphicx, color}

\newtheorem{theorem}{Theorem}[section]
\newtheorem{proposition}[theorem]{Proposition}

\theoremstyle{definition}

\theoremstyle{remark}
\newtheorem{remark}[theorem]{Remark}

\numberwithin{equation}{section}

\newcommand{\CC}{\boldsymbol{C}}
\newcommand{\RR}{\boldsymbol{R}}
\newcommand{\ZZ}{\boldsymbol{Z}}

\newcommand{\lra}{\longrightarrow}
\newcommand{\DS}{\displaystyle}
\newcommand{\op}{\mathrm}

\newcommand{\pa}{\partial}

\newcommand{\ol}{\overline}
\newcommand{\ul}{\underline}

\newcommand{\ali}{\arraycolsep1pt\begin{array}{rl}}
\newcommand{\eali}{\end{array}}
\begin{document}

\title{On origami embeddings of flat tori}

\author{Takashi TSUBOI}
\address{Musashino Center of Mathematical Engineering, Musashino University / RIKEN Interdisciplinary Theoretical and Mathematical Sciences Program
}
\email{mail@tsuboi-takashi.sakura.ne.jp}
\thanks{This work is supported by Musashino University Research Grants.}



\begin{abstract}
We give explicit origami embeddings of a 2-dimensional flat torus of any modulus
in the 3-dimensional Euclidean space.
\end{abstract}

\maketitle

\section{Introduction}
When one learns first courses of Riemannian geometry and the definition of flat tori, one is certain that flat tori do not have $C^2$ isometric embeddings in the 3-dimensional Euclidean space. This can be shown by finding an orthogonal projection to a line which is a Morse function, and then seeing that the points of maxima and minima should have positive Gaussian curvature. Then after having shocked by learning the famous Nash-Kuiper theorem (\cite{nash}, \cite{kuiper}) which asserts any 2-dimensional orientable Riemannian manifold can be $C^1$ isomatrically embedded in the 3-dimensional Euclidean space, and watching the video (\cite{BJLRT}) by V. Borrelli, S. Jabrane, F. Lazarus, D. Rohmer and B. Thibert (see also \cite{BJLT}, \cite{BJLT2}, \cite{EM}), one naturally asks how about origami embeddings. Now, a lot of origami embeddings of flat tori are known. In the litterature, one can find Zalgaller \cite{Zalgaller}, Bern-Hayes \cite{BH}, Segerman \cite{Segerman_book}, etc. There are also those found on the webpages as \cite{Lamb}, \cite{Segerman_video}, \cite{Segerman_shop}, \cite{Tachi}, etc. which attracted the author to think about the moduli of origami embedded flat tori. 

Thus the question treated in this paper is whether a 2-dimensional flat torus of any modulus can be origami embedded in the 3-dimensional Euclidean space. 
The moduli space of flat tori can be obtained from its fundamental domain which is $\{z\in \CC\ \big| \ |\op{Re}\,z|\leqq \DS\frac{1}{2},\ |z|\geqq1\}$. Zalgaller \cite{Zalgaller} explicitly gave origami embeddings of flat tori with moduli of large $\op{Im}\,z$.
Thus the following theorem may be known for people who know the papers cited above. 
However the author thinks that there are still several interesting mathematical feature in the proof given in this article.

\begin{theorem}\label{th:main}
A 2-dimensional flat torus of any modulus can be origami embedded in the 3-dimensional Euclidean space.
\end{theorem}

We prove Theorem \ref{th:main} in the following way.
In Section \ref{sect:annuli}, we give origami embedded flat annuli whose boundaries are the boundaries of regular polygons on two parallel planes with the line joining their centers being orthogonal to the planes. This construction for an equilateral triangle is given by Zalgaller in the section 6 of \cite{Zalgaller}, though the author thinks that it should have been known to craftsmen in old days. Note as in \cite{Zalgaller} that this origami embedding can be deformed by rotating one regular polygon with respect to the other around the perpendicular axis.
Then by putting together two such origami embedded anuuli whose boundaries coincide we obtain an origami immersed flat torus. In Section \ref{sect:tori}, we see the conditions for the two origami embedded anuuli to give rise to an origami embedding. These origami embedded flat tori were explicitly given by Segerman (\cite{Segerman_video}, \cite{Segerman_book}, \cite{Segerman_shop}). In Section \ref{sect:moduli_of_tori}, we calculate the moduli of origami embedded flat tori constructed in this way. We find that the moduli space of such flat tori is related to the space of tangents of cycloids. There we see that this construction gives origami embeddings of flat tori except those with pure imaginary moduli, i.e., except those with rectangular fundamental domains. 
In Section  \ref{sect:rectangular}, we use half of the origami embedded flat tori in Section  \ref{sect:tori} and take the double of them to construct origami embeddings of flat tori with rectangular fundamental domains. This doubling construction is essentially the same as the bending construction for triangular cylinders given by Zalgaller \cite{Zalgaller}.

The rough idea of constructions of origami 
embeddings of flat tori in this article was explained in Japanese in \cite{Tsuboi} where we showed typical examples which are included also in this article.

The author thanks Musashino Center of Mathematical Engineering of Musashino University and RIKEN Interdisciplinary Theoretical and Mathematical Sciences Program for their support. He is grateful to the valuable discussion with Kazuo Masuda and with Shizuo Kaji during the preparation of this article.

\begin{figure}
\begin{center}
\includegraphics[width=3cm]{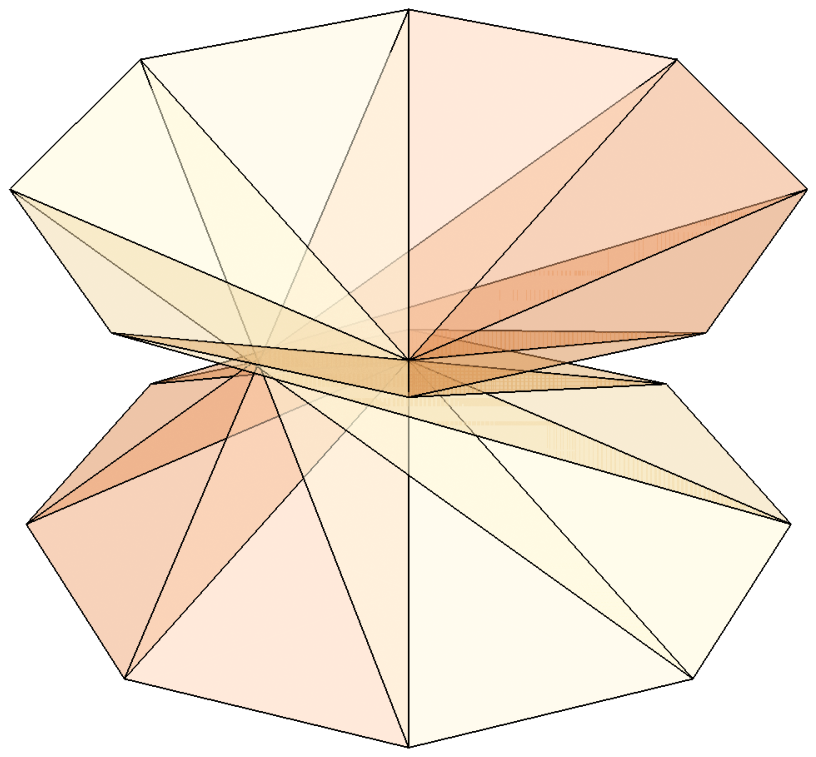}
\includegraphics[width=3cm]{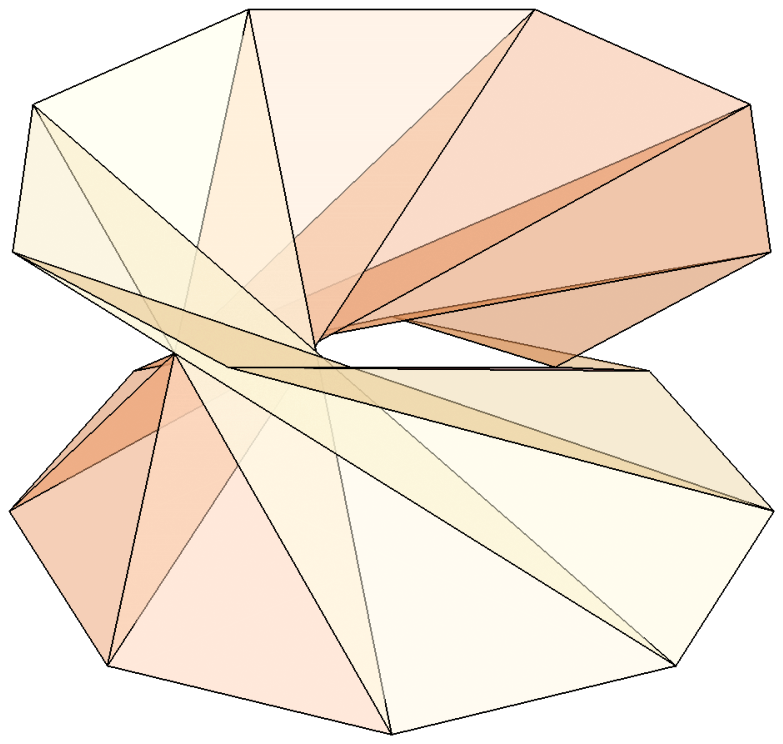}
\includegraphics[width=3cm]{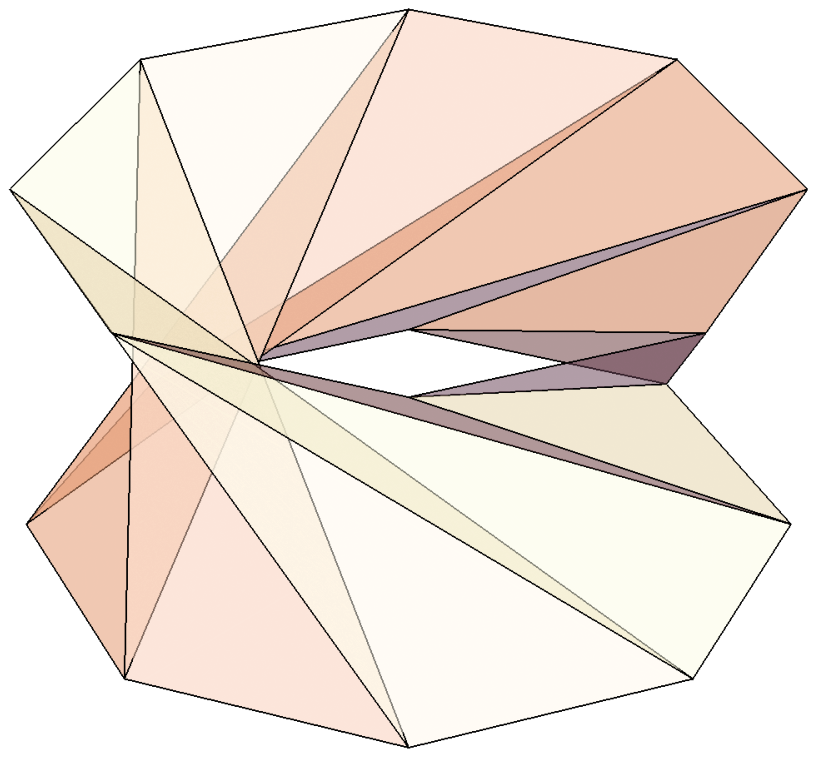}

\includegraphics[width=3cm]{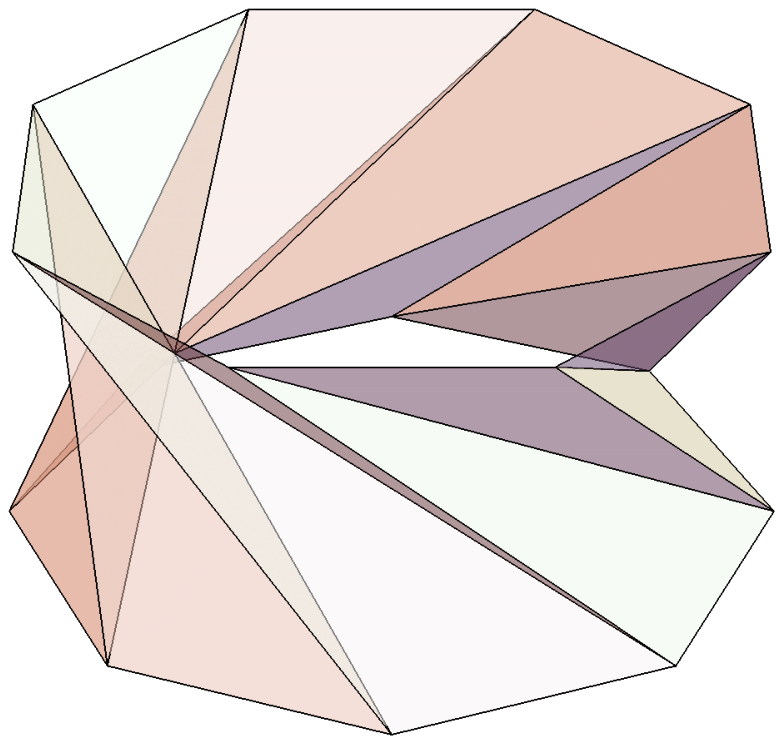}
\includegraphics[width=3cm]{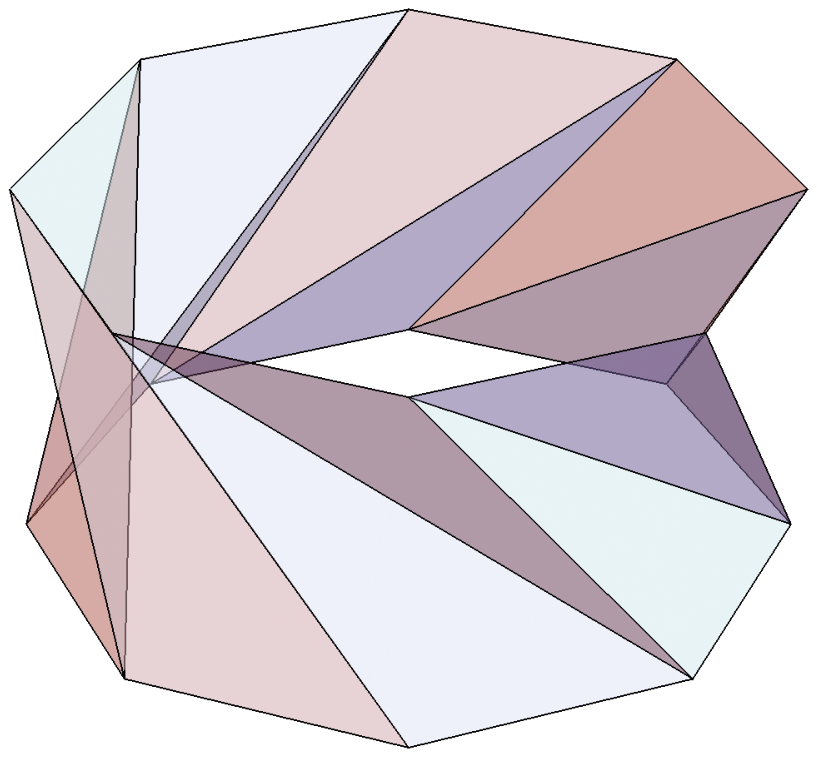}
\includegraphics[width=3cm]{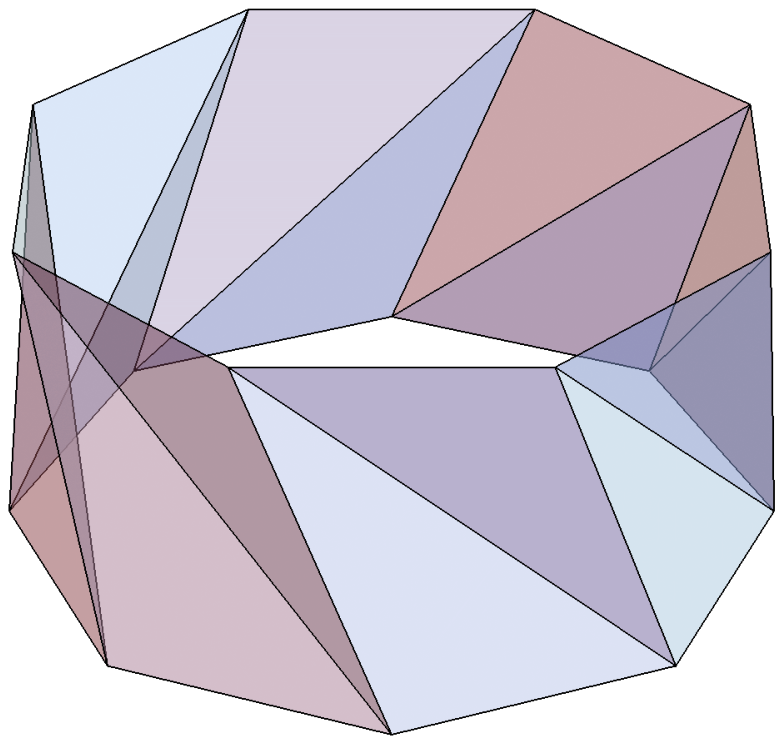}

\includegraphics[width=3cm]{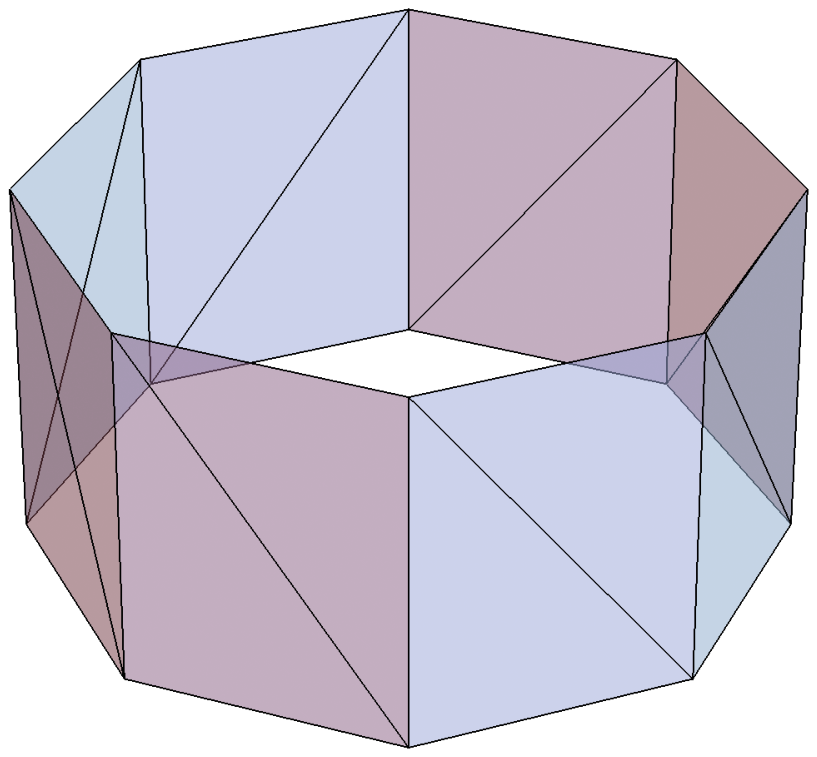}
\includegraphics[width=3cm]{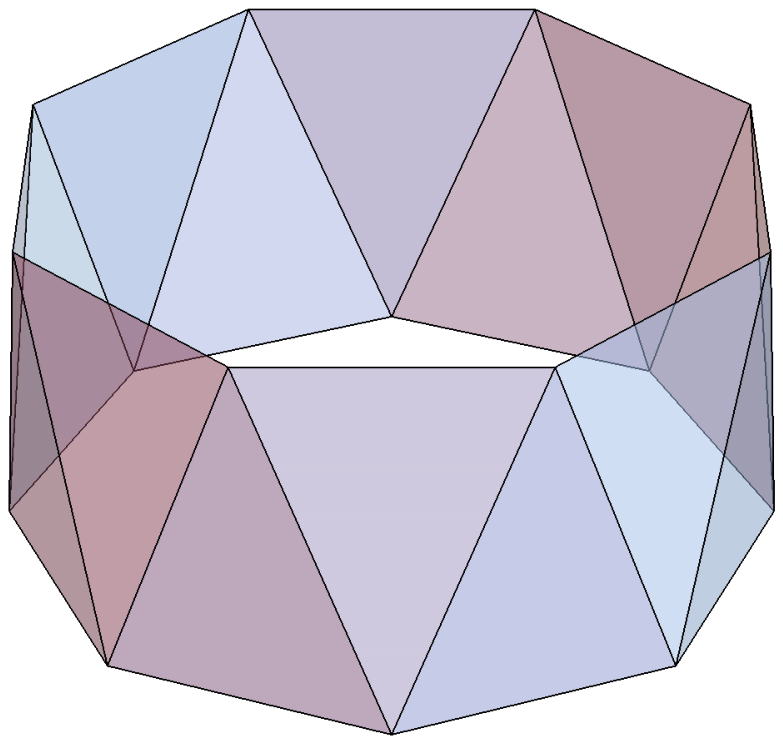}
\includegraphics[width=3cm]{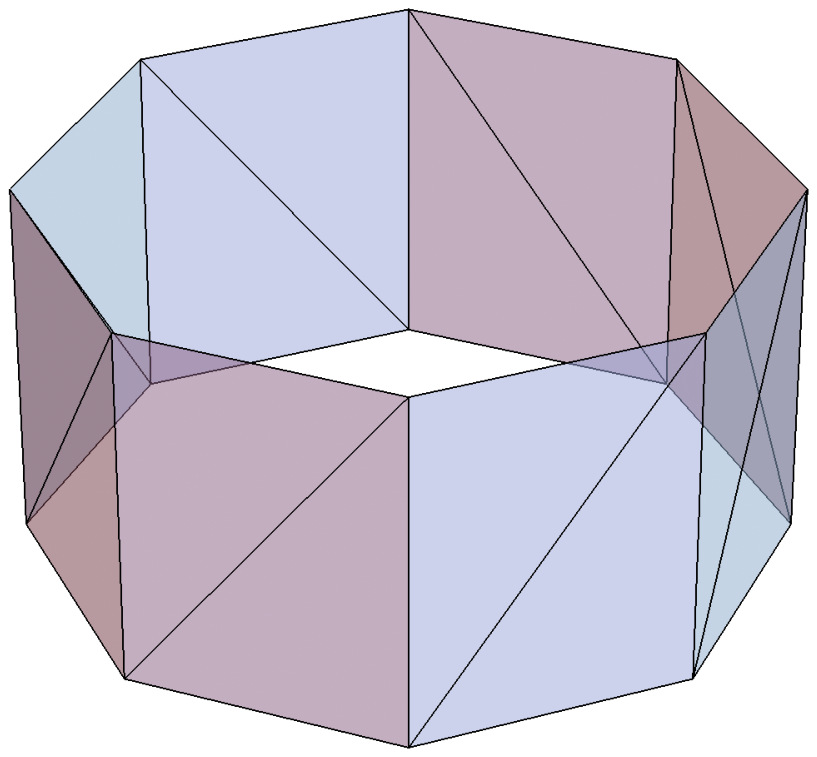}

\includegraphics[width=3cm]{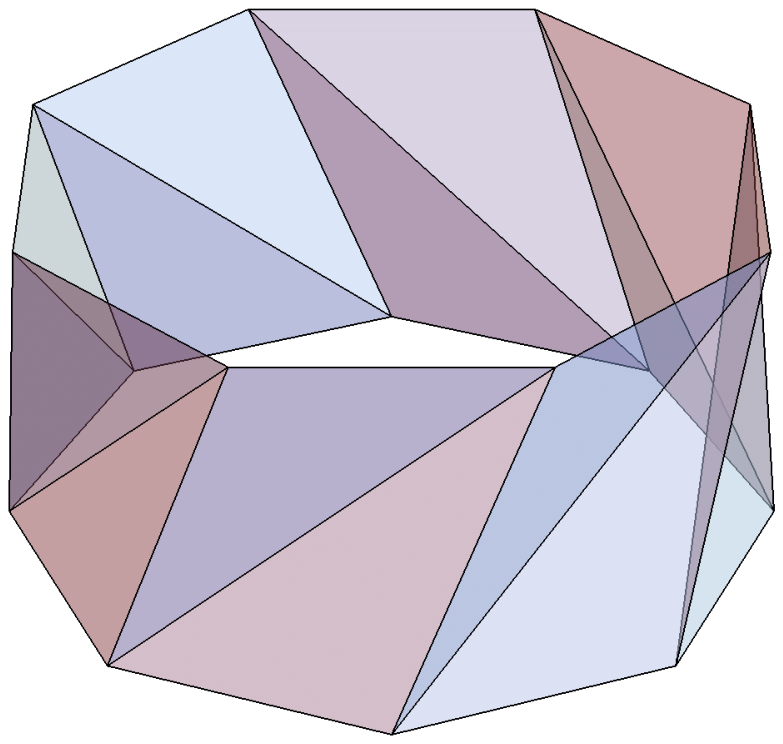}
\includegraphics[width=3cm]{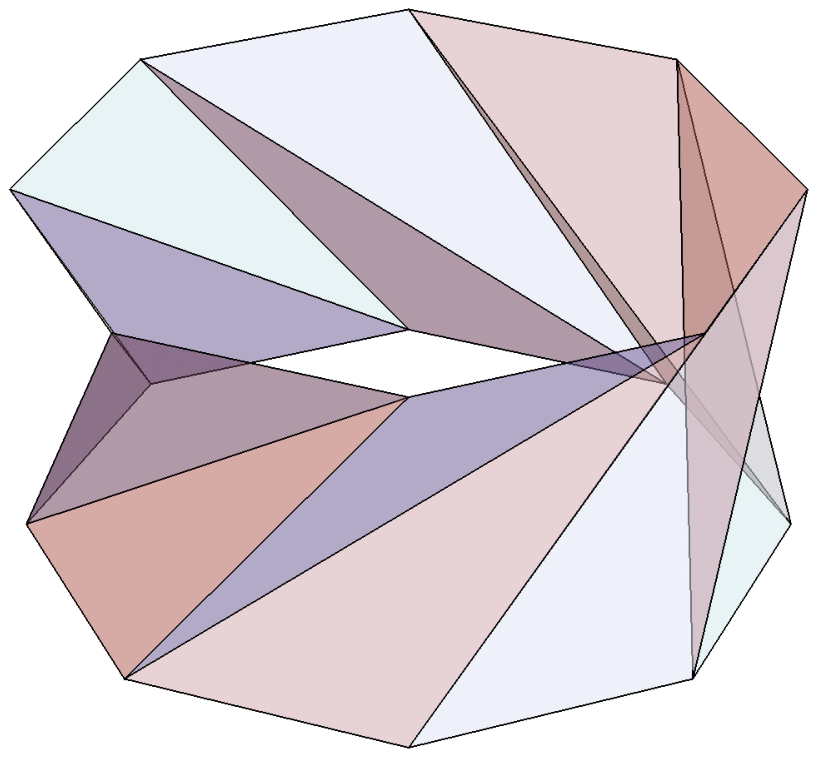}
\includegraphics[width=3cm]{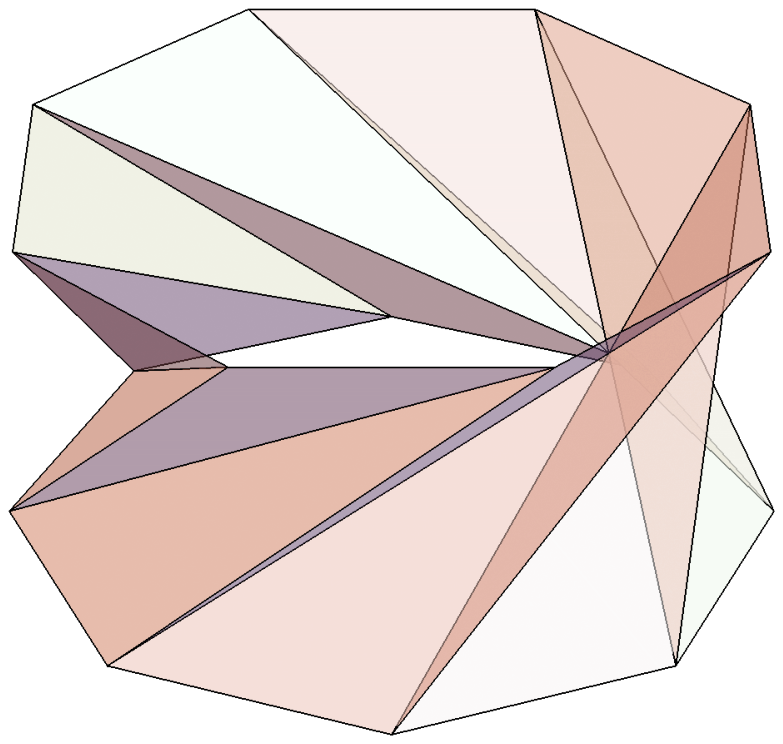}

\includegraphics[width=3cm]{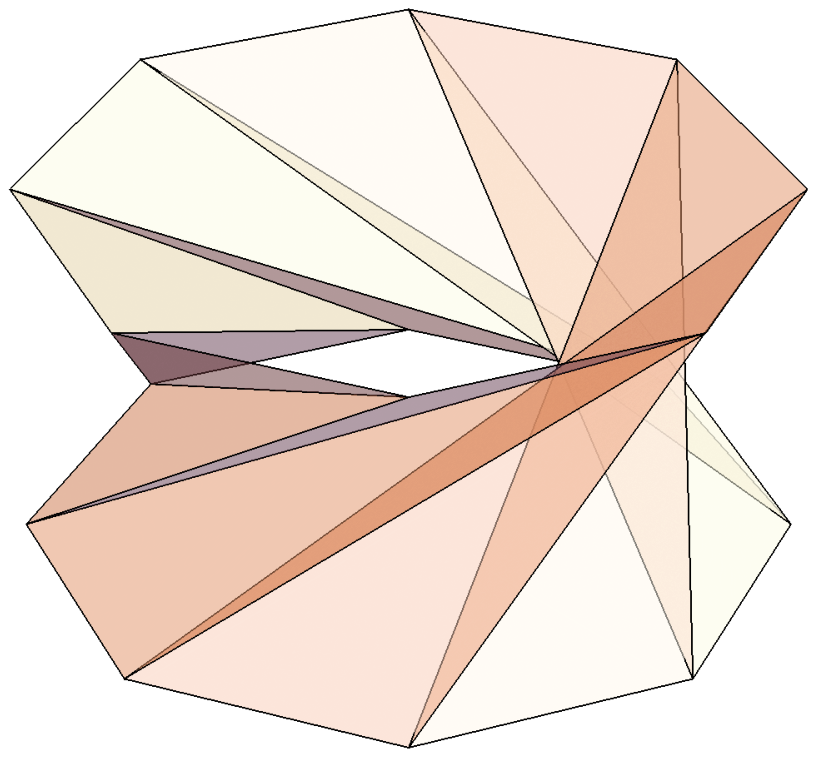}
\includegraphics[width=3cm]{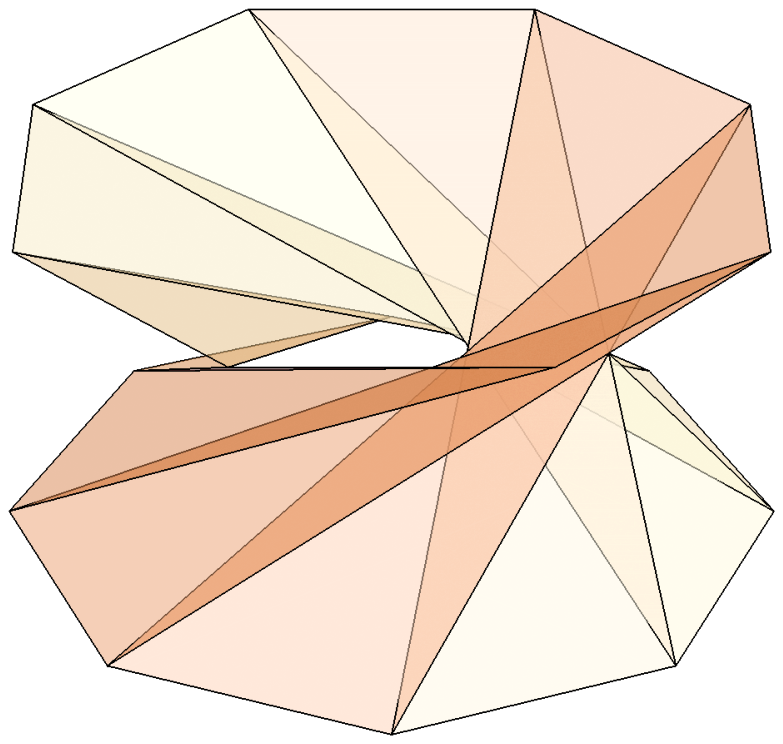}
\includegraphics[width=3cm]{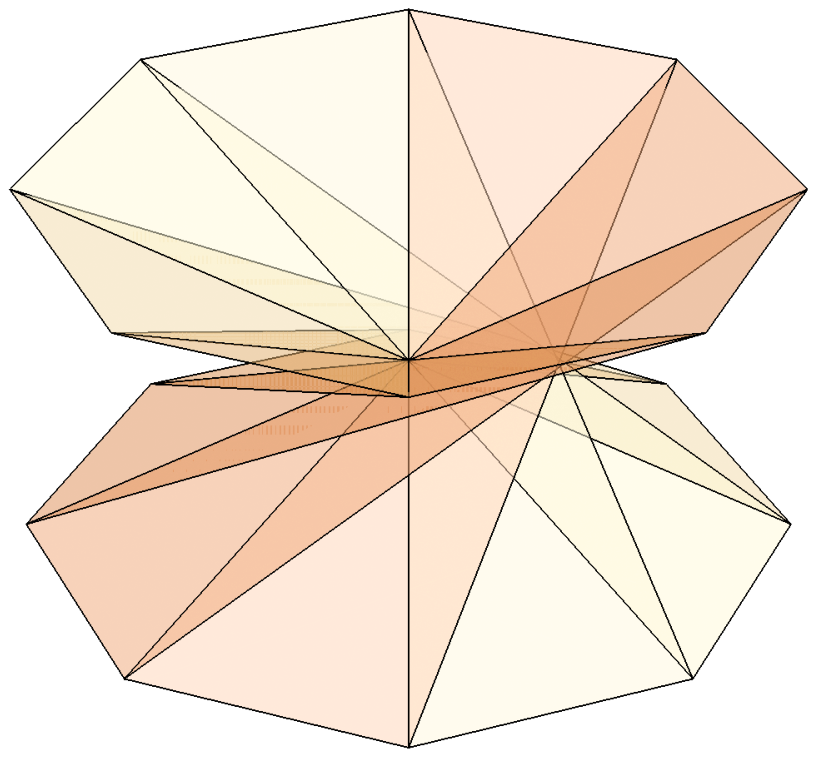}
\end{center}
\caption{\small $A_n^\rho$ for $h=1$, $n=8$ and $\rho =
\DS\frac{k}{2n}$ ($k=-n$, $-n+1$, \dots, $n-2$).
The left figure of the uppermost line shows $A_8^{-8/16}$ and the right figure of the lowermost line shows 
$A_8^{+6/16}$, and these two are not embedded annuli. }\label{fig:A_8_rho}
\end{figure}

\begin{figure}
\begin{center}
\includegraphics[width=4cm]{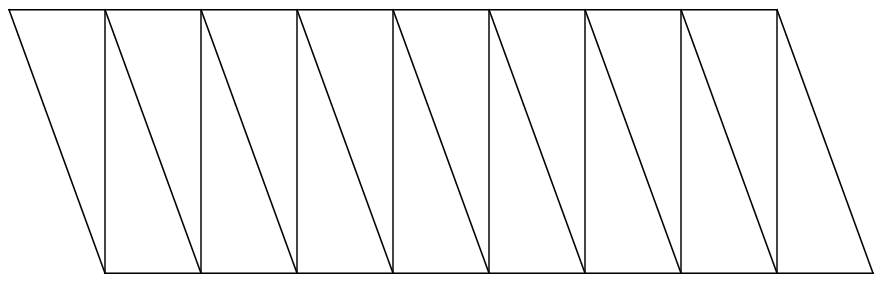}
\includegraphics[width=4cm]{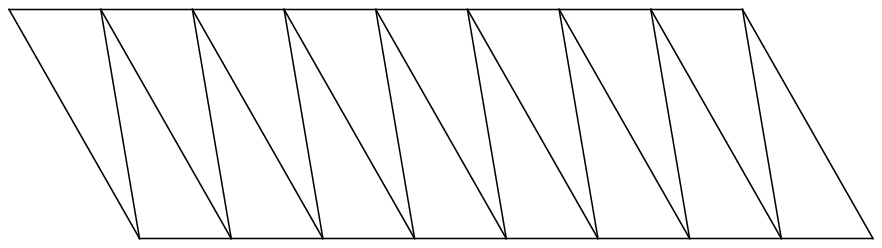}
\includegraphics[width=4cm]{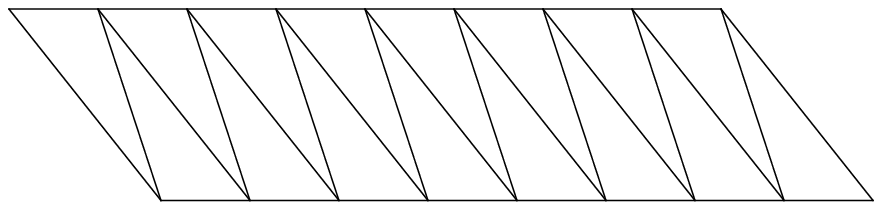}
\vskip4mm

\includegraphics[width=4cm]{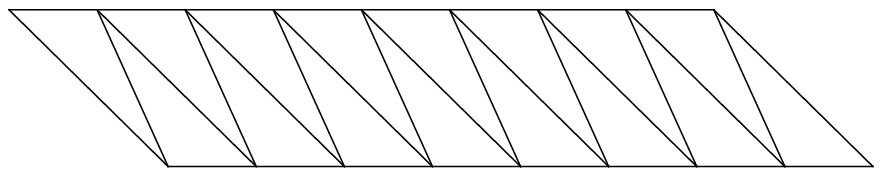}
\includegraphics[width=4cm]{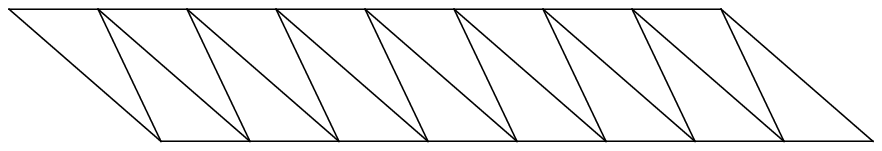}
\includegraphics[width=4cm]{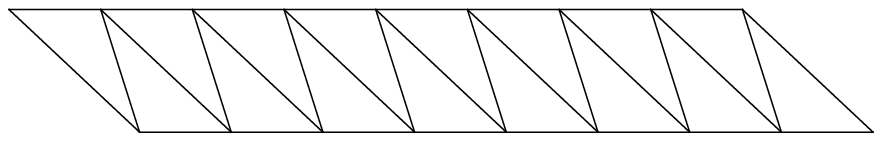}
\vskip4mm

\includegraphics[width=4cm]{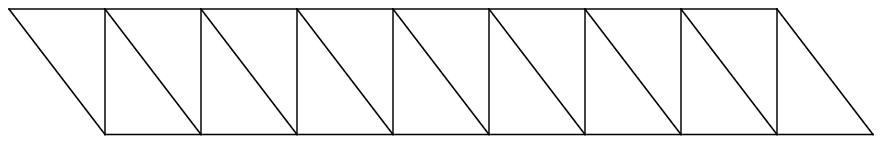}
\includegraphics[width=4cm]{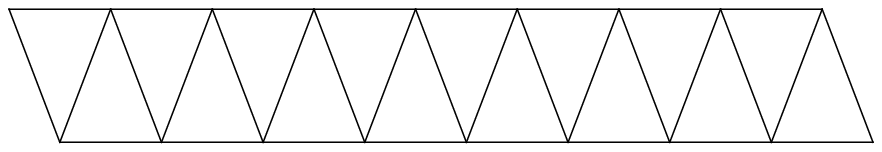}
\includegraphics[width=4cm]{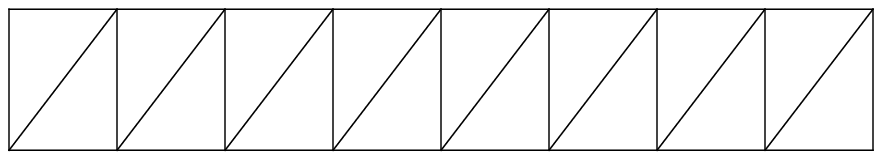}
\vskip4mm

\includegraphics[width=4cm]{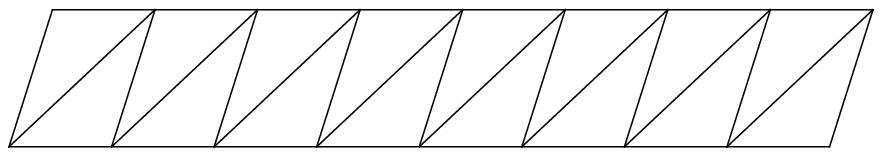}
\includegraphics[width=4cm]{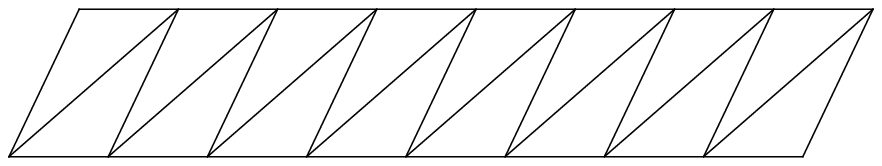}
\includegraphics[width=4cm]{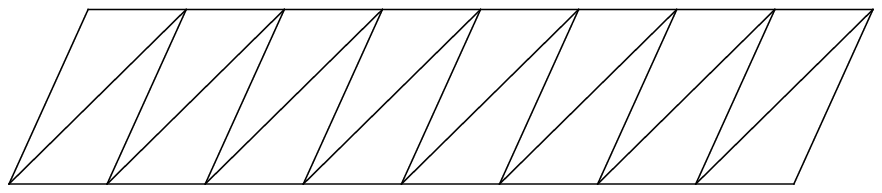}
\vskip4mm

\includegraphics[width=4cm]{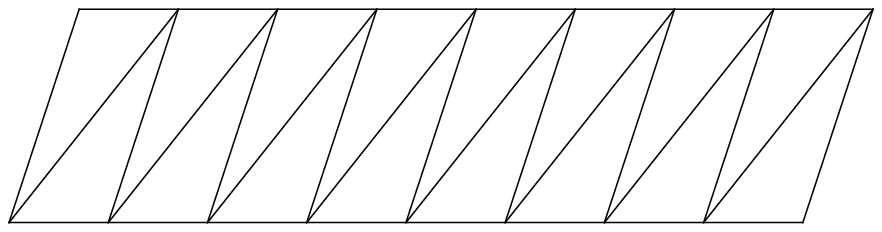}
\includegraphics[width=4cm]{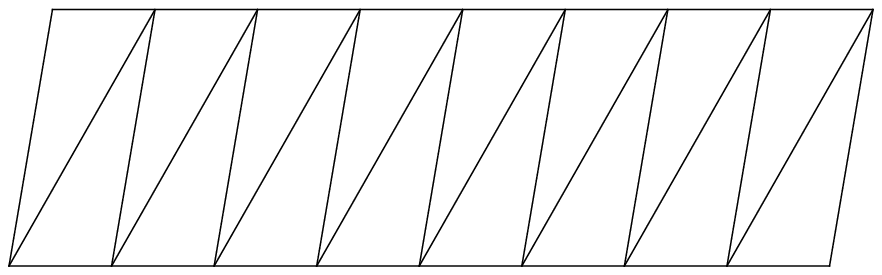}
\includegraphics[width=4cm]{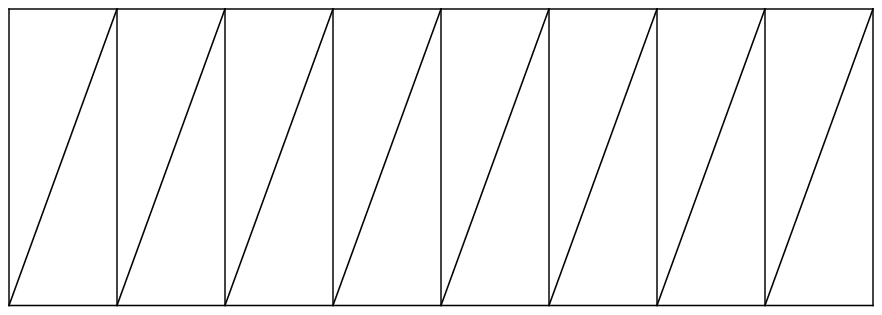}
\end{center}
\caption{\small The developments of $A_8^\rho$ in Figure \ref{fig:A_8_rho} ($h=1$, $n=8$, $\rho=\DS\frac{k}{2n}$; $k=-n$, $-n+1$, \dots, $n-2$). The development for $\rho=\DS\frac{-8}{16}$ shown at the left figure of the uppermost line and that for $\rho=\DS\frac{6}{16}$ shown at the right figure of the lowermost line do not correspond to embeddings.}\label{fig:DevA_8_rho}
\end{figure}

\section{Construction of embeddings of flat annuli}\label{sect:annuli}
In this section we construct origami embedded flat annuli.

Let us use the coodinates $\CC\times\RR$ for the 3-dimensional
Euclidean space.
Let $n$ be an integer greater than 2.
First, consider the regular $n$-gon $P_0\cdots P_{n-1}P_n$ ($P_n=P_0$) on $\CC\times\{0\}\subset \CC\times\RR$ with vertices $P_k=(e^{2\pi \sqrt{-1}k/n},0)$ ($k=0$, \dots, $n$).
Secondly, for a positive real number $h$ and an element $\rho\in \RR/\ZZ$,
consider the regular $n$-gon $Q_0^\rho\cdots Q_{n-1}^\rho Q_n^\rho$ ($Q_n^\rho=Q_0^\rho$) on 
$\CC\times\{h\}\subset \CC\times\RR$ with vertices 
$Q_k^\rho=(e^{2\pi\sqrt{-1}\rho+2\pi\sqrt{-1}k/n},h)$ ($k=0$, \dots, $n$).

Consider the triangles $\triangle P_0P_1Q_1^\rho$ and $\triangle Q_{0}^\rho Q_{1}^\rho P_{1}$. Hereafter we take a representative $\rho\in \RR$ of $\rho \in \RR/\ZZ$ and we assume that $$-\DS\frac{1}{2} <\rho  < \DS\frac{1}{2} -\DS\frac{1}{n}.$$
Note first that the lengths of the edges $P_0P_1$ and $Q_0^\rho Q_1^\rho$ are equal as the edges of the congruent regular $n$-gons. Secondly, the lengths of the edges $P_0Q_0^\rho$ and $P_1Q_{1}^\rho$ are equal because they are mapped by the rotation by $\DS\frac{2\pi}{n}$ around the real axis $\{0\}\times \RR$. Since the triangles $\triangle P_0P_1Q_1^\rho$ and $\triangle Q_{0}^\rho Q_{1}^\rho P_{1}$ share the edge $P_1Q_1^\rho$, they are congruent; 
$\triangle P_0P_1Q_1{}^u\cong\triangle Q_1{}^uQ_0{}^uP_0$. 
We rotate these triangles around the real axis $\{0\}\times \RR$ by $\DS\frac{2k\pi}{n}$ and obtain $2n$ triangles $\triangle P_kP_{k+1}Q_{k+1}^\rho$, $\triangle Q_{k+1}^\rho Q_k^\rho P_{k}$ ($k=0$, \dots $n-1$). Let $A_n^\rho$ denote the union of them;
$$A_n^\rho=\bigcup_{k=0,\dots,n-1}\triangle P_kP_{k+1}Q_{k+1}^\rho\cup\triangle Q_{k+1}^\rho Q_k^\rho P_{k}.$$
It is necessary to verify that the interiors of these triangles are disjoint.
We find the infimum and the supremum of $\rho$ 
such that the interiors of these triangles are disjoint, and they are $-\DS\frac{1}{2}$ and 
$\DS\frac{1}{2}-\DS\frac{1}{n}$, respectively. 
At $\rho=-\DS\frac{1}{2}$ the edges $P_kQ_k^\rho$ ($k=0$, \dots $n-1$) 
pass through the center 
$(0,\DS\frac{h}{2})$ of $A_n^\rho$, 
and at $\rho=\DS\frac{1}{2}-\DS\frac{1}{n}$ so do the edges 
$P_kQ_{k+1}^\rho$ ($k=0$, \dots $n-1$).
Thus the interiors of these triangles are disjoint by the assumption
that $-\DS\frac{1}{2} <\rho  < \DS\frac{1}{2} -\DS\frac{1}{n}$.

If $\rho\in[-\DS\frac{1}{n},0]$,
$A_n^\rho$ is the union of side faces of the convex hull of the union of the two regular $n$-gons $P_0\cdots P_{n-1}P_n$ and $Q^\rho_0\cdots Q^\rho_{n-1}Q^\rho_n$. If $\rho$ is $-\DS\frac{1}{n}$ or $0$, $A_n^\rho$ is the union of the side faces of the regular $n$ angular prism and if $\rho= -\DS\frac{1}{2n}$, $A_n^\rho$ is the union of side faces of the uniform $n$-gonal antiprism. See Figure \ref{fig:A_8_rho}.

If we develope $A_n^\rho$, we obtain a parallelogram. For, since $\triangle P_0P_1Q_1^\rho\cong\triangle Q_1^\rho Q_0^\rho P_0$, 
$\triangle P_0P_1Q^\rho_1\cup\triangle Q^\rho_1 Q^\rho_0P_1$ is a parallelogram. As a union of $n$ copies of this parallelogram,
$A_n^\rho$ is developed to a parallelogram 
with the edge of length $2n\sin{\DS\frac{\pi}{n}}$ and height equals to the height of the triangle 
$\triangle P_0P_1Q_1^\rho$ with respect to
the base edge $P_0P_1$. Thus $A_n^\rho$ is an origami emmbeded flat annulus.
See Figure \ref{fig:DevA_8_rho}.

Since the length of edges of $\triangle P_0P_1Q_1^\rho$ are computed as $P_0P_1=2\sin(\DS\frac{\pi}{n})$,  
$P_1Q_1^\rho=\sqrt{4(\sin(\pi\rho))^2+h^2}$ and $P_0Q_1^\rho=\sqrt{4(\sin(\pi\rho+\DS\frac{\pi}{n}))^2+h^2}$,
the point $X_1^\rho$ on the line $P_0P_1$ which is $Q_1^\rho$ projected to the line $P_0P_1$ satisfies that
$$\ali
{P_0X_1^\rho}=&\DS\frac{(P_0P_1)^2+(P_0Q_1^\rho)^2-(P_1Q_1^\rho)^2}{2P_0P_1}
\\=&2\cos(\pi\rho)\sin(\pi\rho+\DS\frac{\pi}{n})=
\sin(2\pi\rho+\DS\frac{\pi}{n})+\sin(\DS\frac{\pi}{n}).
\eali$$


Then the height of the triangle $\triangle P_0P_1Q^\rho_1$ with respect to the base edge $P_0P_1$ is equal to
$$
\ali X_1^\rho Q_1^\rho=&\sqrt{(P_0Q_1^\rho)^2
-(P_0X_1^\rho)^2}
\\=&\sqrt{h^2+4(\sin(\pi\rho))^2(\sin(\pi\rho+\DS\frac{\pi}{n}))^2}
=\sqrt{h^2+\big(\cos(2\pi\rho+\DS\frac{\pi}{n})-\cos\DS\frac{\pi}{n}\big)^2}.

\eali
$$

\begin{figure}\begin{center}
\includegraphics[height=80mm]{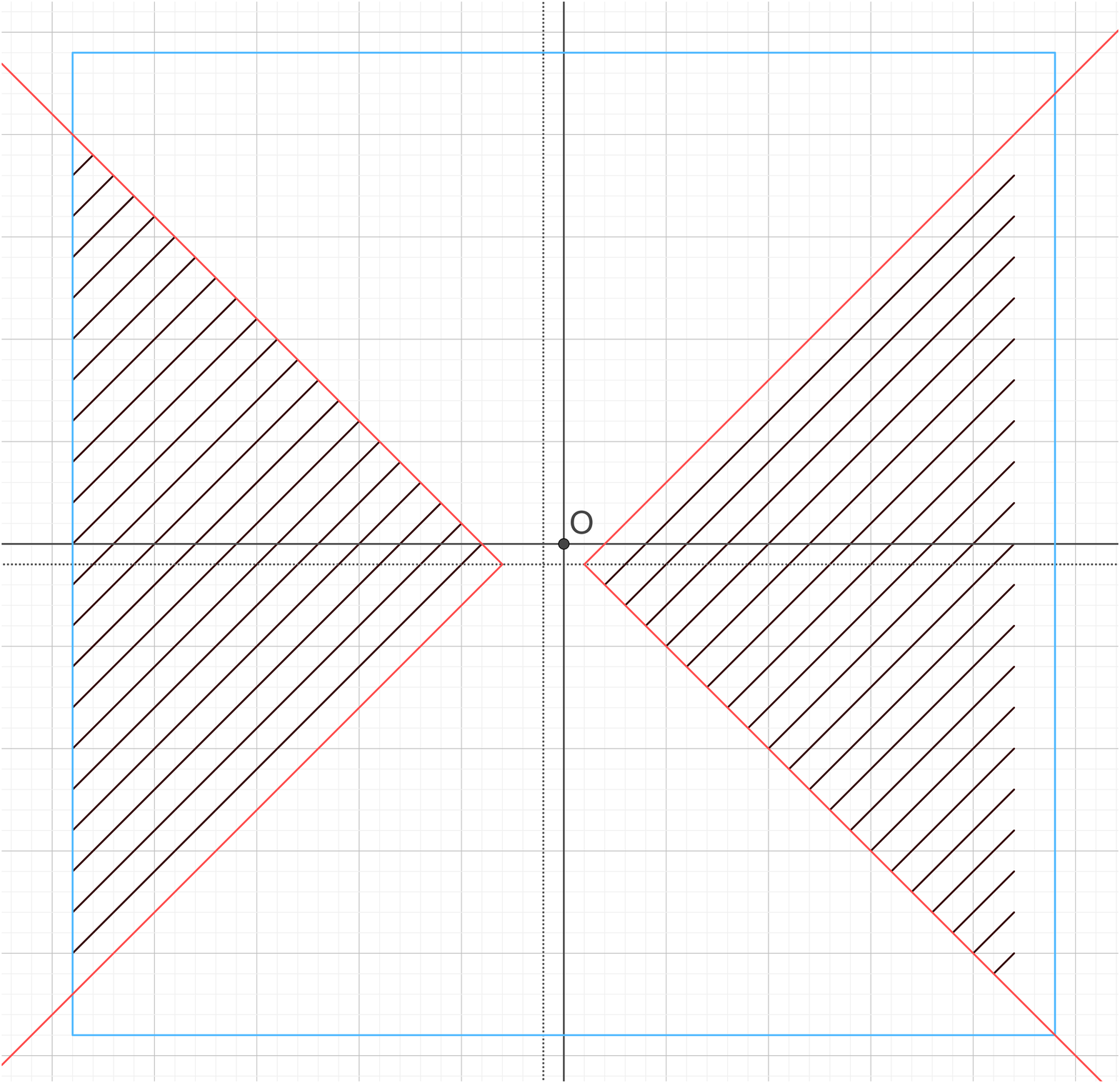}
\caption{\small
For $A_n^\rho$ and $A_n^\sigma$ with
$\sigma=\rho+\DS\frac{\ell}{n}$, the condition that 
 $A_n^\rho$ is closer to $\{0\}\times \RR$ than $A_n^\sigma$ is in the interior of the black segments of slope 1 in the $\rho\sigma$-coordinates. 
$n=24$ and the blue square is $[-\DS\frac{1}{2},\DS\frac{1}{2}]\times[-\DS\frac{1}{2},\DS\frac{1}{2}]$.
}\label{fig:rhosigma}　　　
\end{center}
\end{figure}

\begin{figure}\begin{center}
\includegraphics[height=60mm]{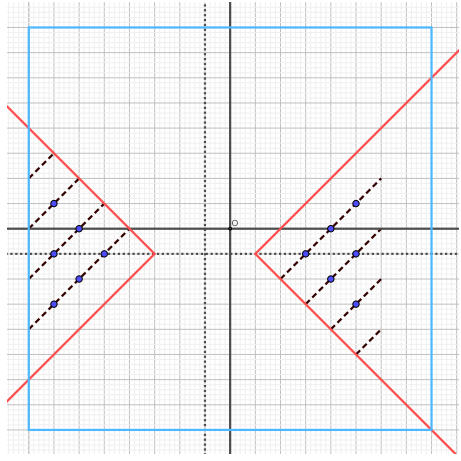}
\caption{\small
The boundaries of 
$A_n^\rho$ ($n=8$, $\rho=\DS\frac{k}{16}$ ($k=-8$, \dots, $6$)) in Figure \ref{fig:A_8_rho} coincide alternately.
The conditions on the boundary ($\sigma=\rho+\DS\frac{\ell}{n}$),
 on the disjointness and on the position to the axis $\{0\}\times\RR$
is satisfied by these 12 points indicated on the $\rho\sigma$-coordinates.
The blue square is $[-\DS\frac{1}{2},\DS\frac{1}{2}]\times[-\DS\frac{1}{2},\DS\frac{1}{2}]$.}
\label{fig:rhosigma8}　　　
\end{center}
\end{figure}

\begin{figure}
\begin{center}
\includegraphics[width=3.5cm]{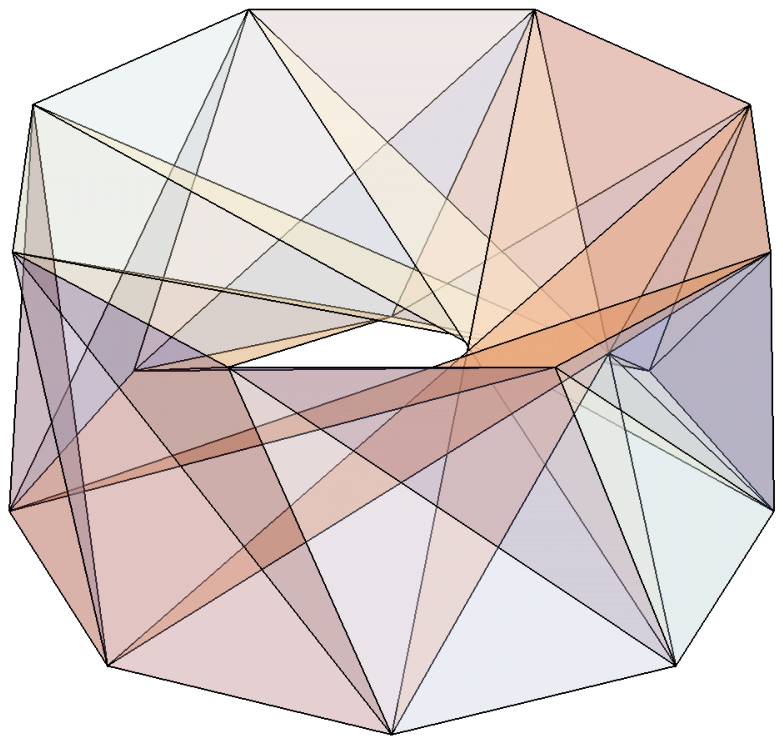}
\includegraphics[width=3.5cm]{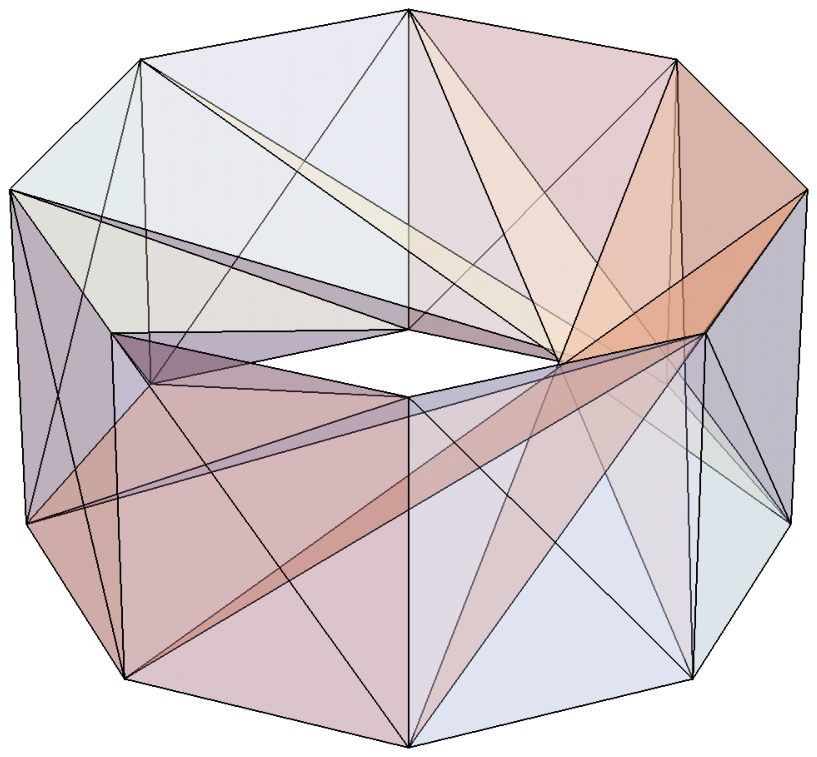}
\includegraphics[width=3.5cm]{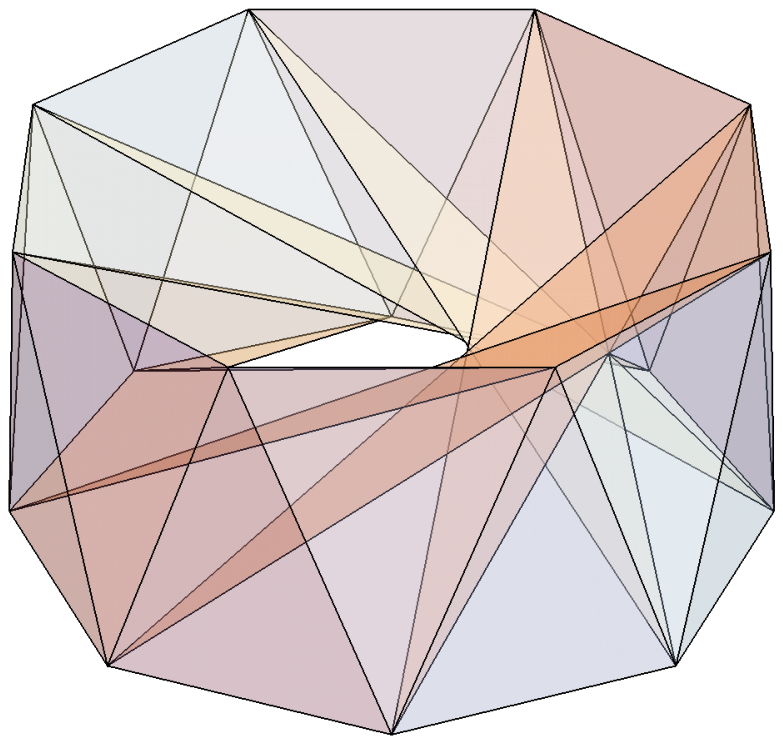}

$(\DS\frac{5}{16},-\frac{3}{16})$\hskip2.4cm 
$\DS(\frac{4}{16}, -\frac{2}{16})$\hskip2cm
$\DS(\frac{5}{16},-\frac{1}{16})$
\vskip3mm

\includegraphics[width=3.5cm]{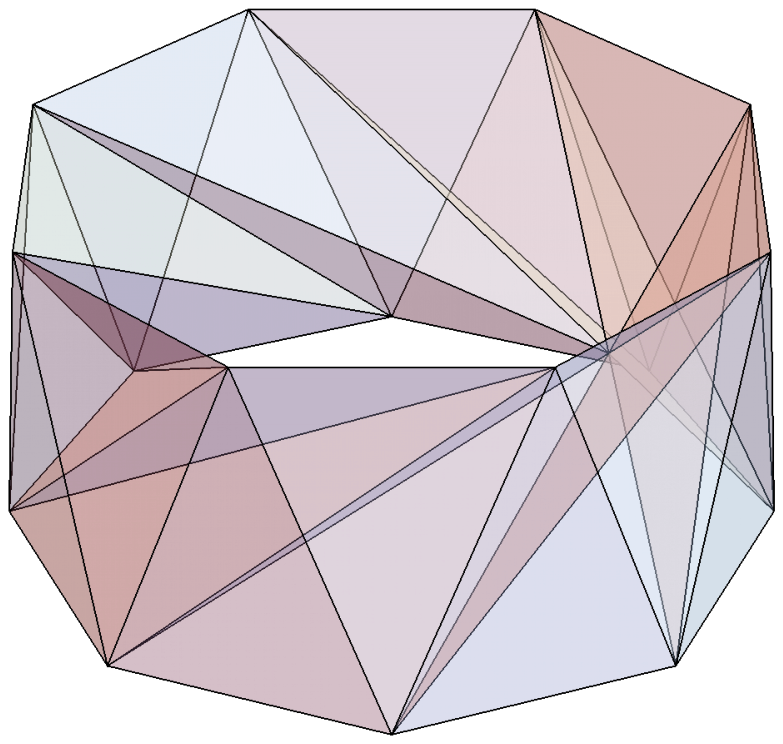}
\includegraphics[width=3.5cm]{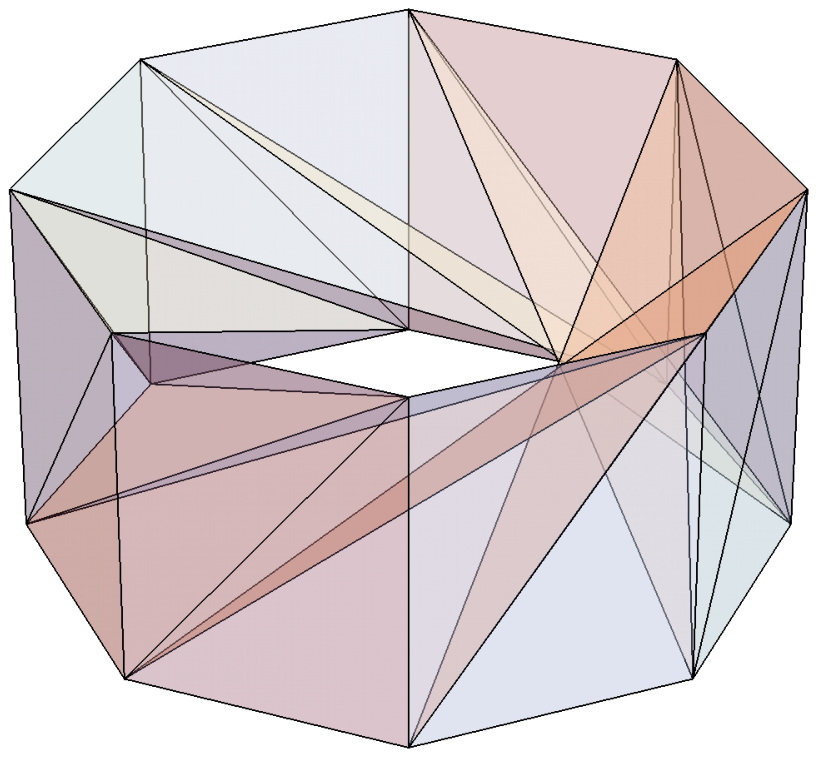}
\includegraphics[width=3.5cm]{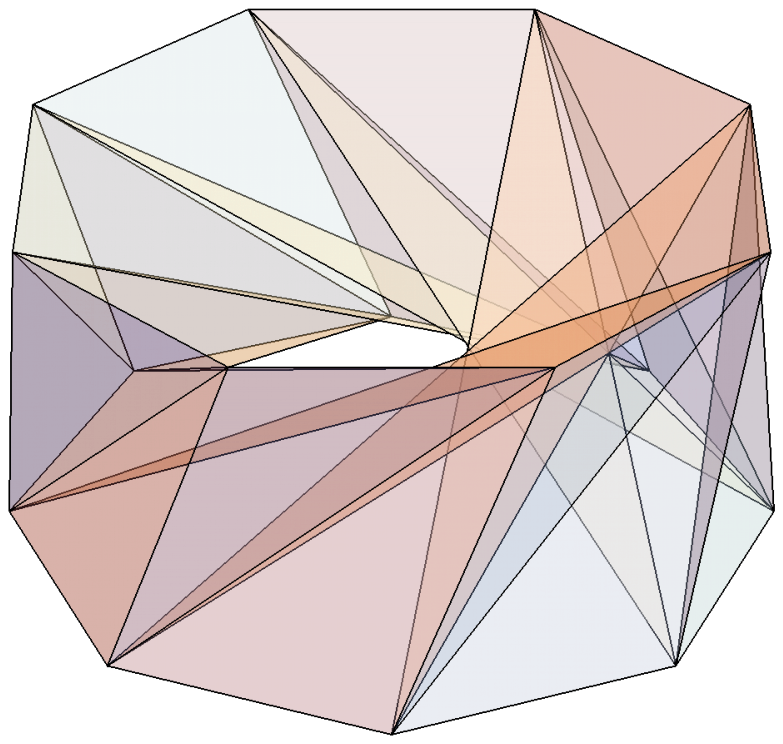}

$\DS(\frac{3}{16},-\frac{1}{16})$\hskip2.4cm
$\DS(\frac{4}{16},0)$\hskip2.4cm
$\DS(\frac{5}{16},\frac{1}{16})$
\vskip3mm

\includegraphics[width=3.5cm]{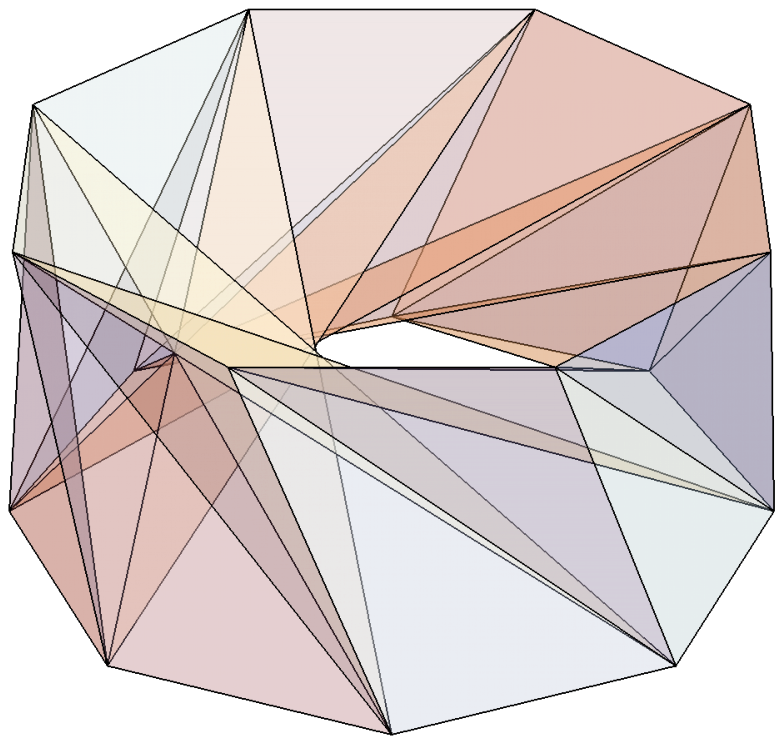}
\includegraphics[width=3.5cm]{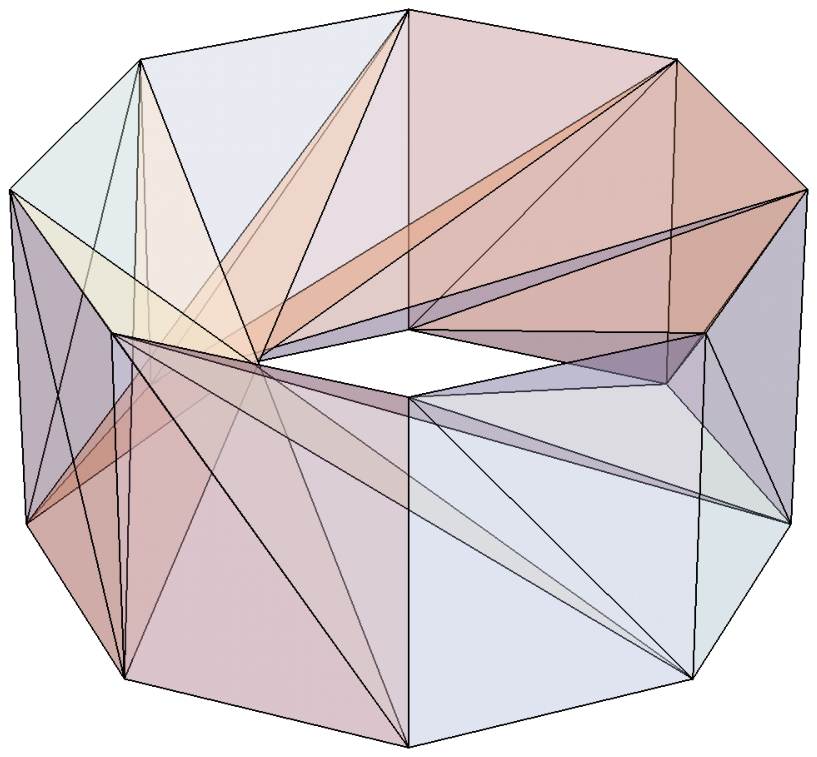}
\includegraphics[width=3.5cm]{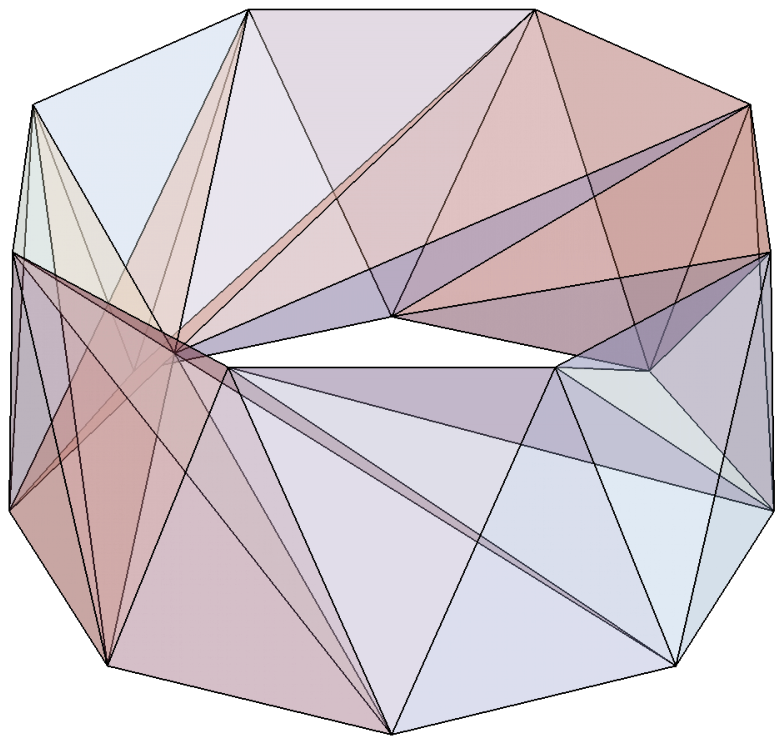}

$\DS(-\frac{7}{16},-\frac{3}{16})$\hskip1.8cm
$\DS(-\frac{6}{16},-\frac{2}{16})$\hskip1.8cm
$\DS(-\frac{5}{16},-\frac{1}{16})$
\vskip3mm

\includegraphics[width=3.5cm]{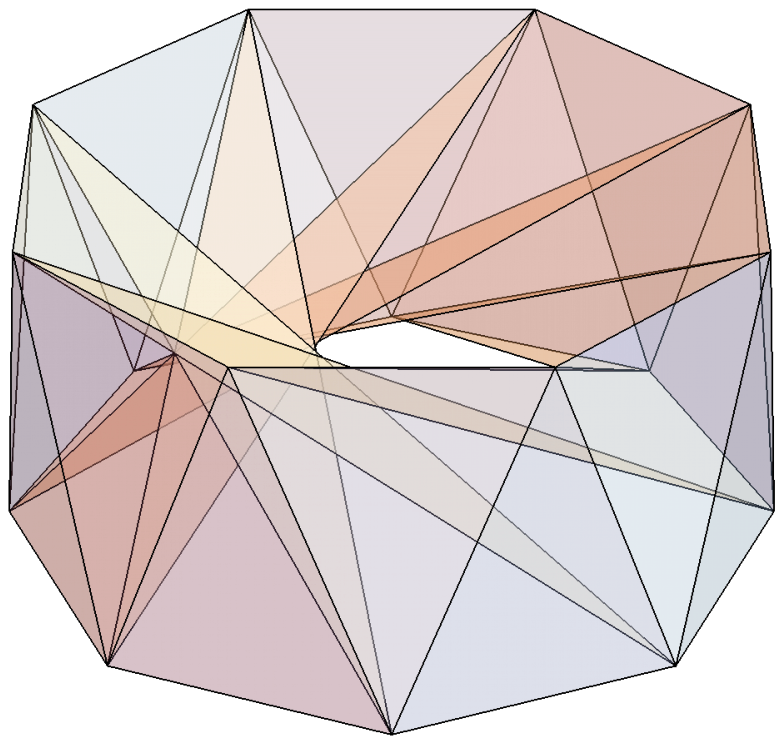}
\includegraphics[width=3.5cm]{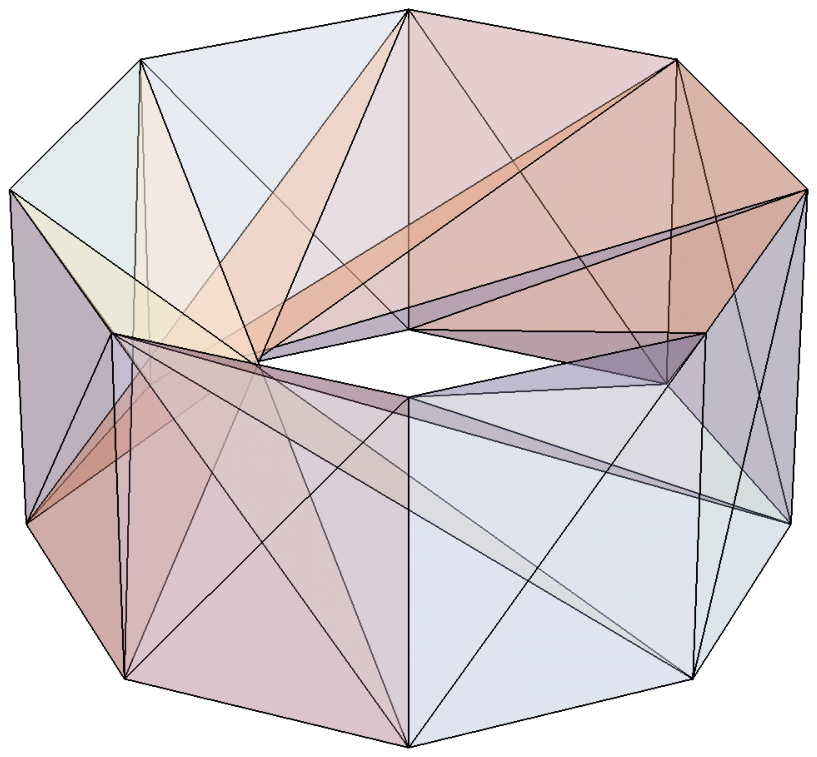}
\includegraphics[width=3.5cm]{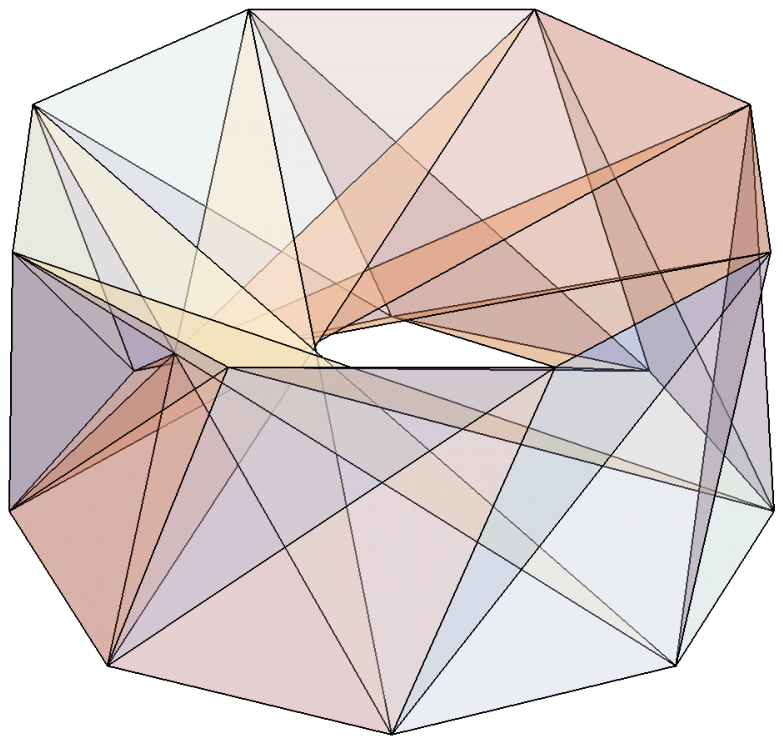}

$\DS(-\frac{7}{16},-\frac{1}{16})$\hskip2cm
$\DS(-\frac{6}{16},0)$\hskip2.1cm
$\DS(-\frac{7}{16},\frac{1}{16})$

\end{center}
\caption{\small
Origami embeddings of flat tori $A_8^\rho\cup A_8^\sigma$ and the coodinates $(\rho,\sigma)$
for 
$n=8$, $h=1$.
}\label{fig:rhosigma8emb}
\end{figure}

\begin{figure}\begin{center}
\includegraphics[height=1.3cm]{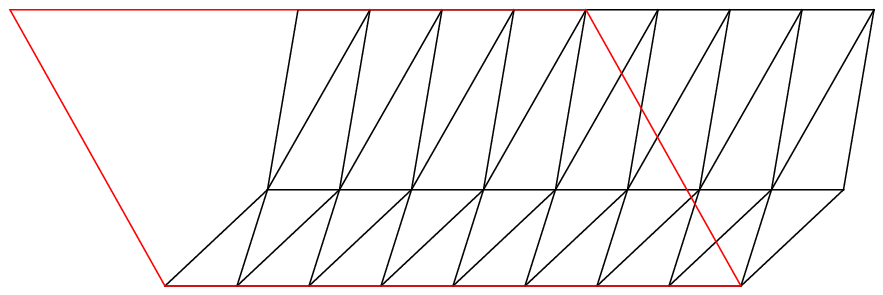}\ \ 
\includegraphics[height=1.3cm]{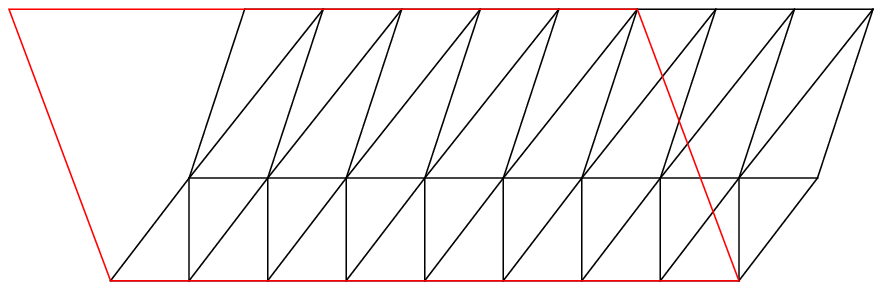}\ \ 
\includegraphics[height=1.3cm]{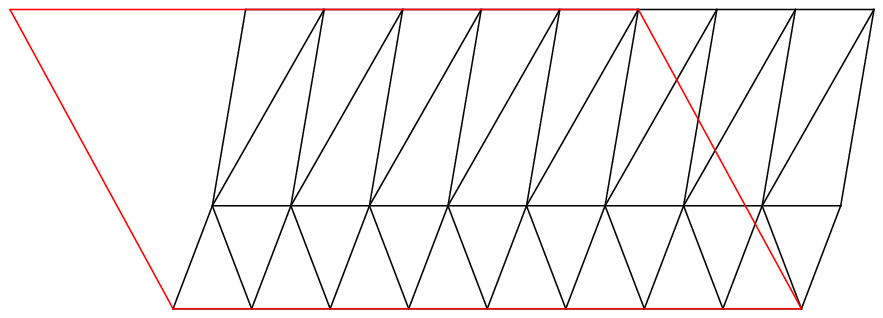}
\vskip1mm

$\DS-\frac{3}{16}-\DS\frac{5}{16}=-\frac{4}{8}$\hskip1.8cm 
$\DS-\frac{2}{16}-\frac{4}{16}=-\frac{3}{8}$\hskip1.6cm
$\DS-\frac{1}{16}-\frac{5}{16}=-\frac{3}{8}$
\vskip3mm

\includegraphics[height=1.3cm]{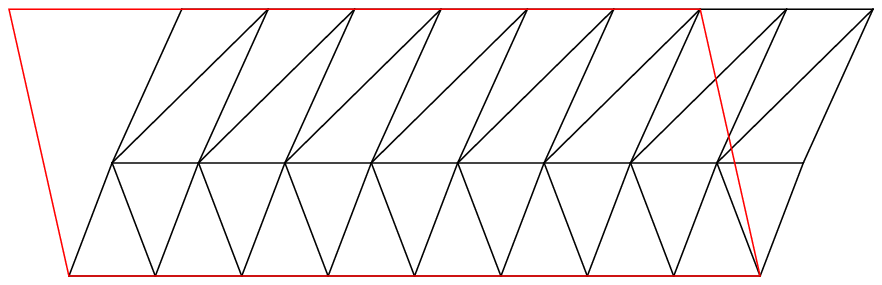}\ \ 
\includegraphics[height=1.3cm]{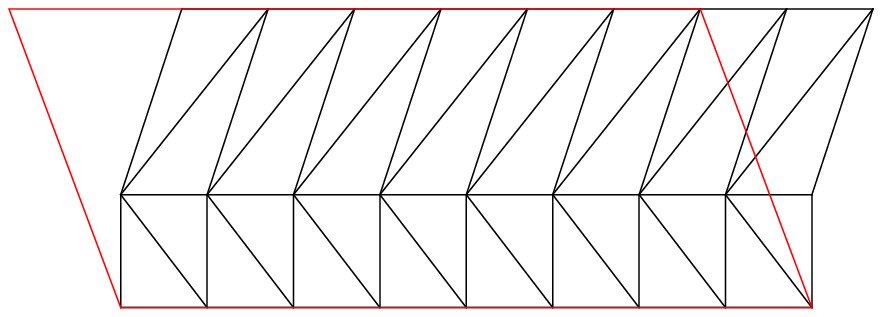}\ \ 
\includegraphics[height=1.3cm]{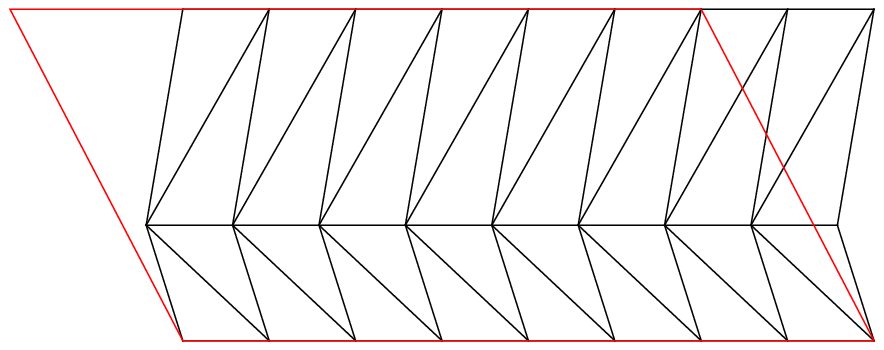}
\vskip1mm

$\DS-\frac{1}{16}-\frac{3}{16}=-\frac{2}{8}$\hskip2.0cm
$\DS0-\frac{4}{16}=-\frac{2}{8}$\hskip1.6cm
$\DS \frac{1}{16}-\frac{5}{16}=-\frac{2}{8}$

\vskip3mm

\includegraphics[height=1.3cm]{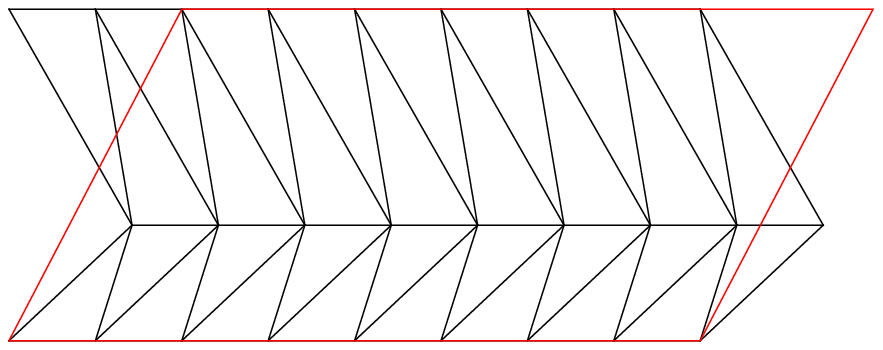}\ \ 
\includegraphics[height=1.3cm]{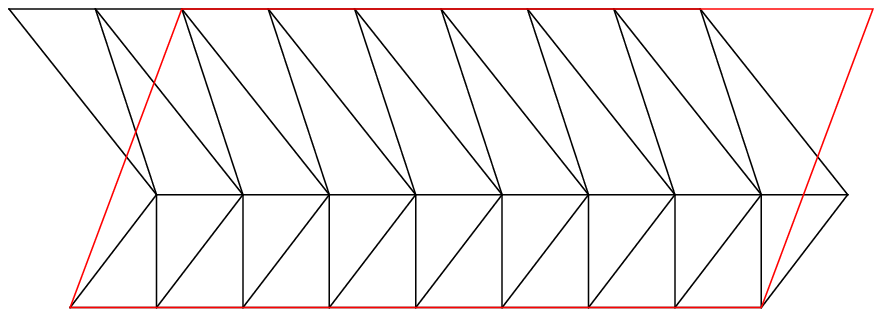}\ \ 
\includegraphics[height=1.3cm]{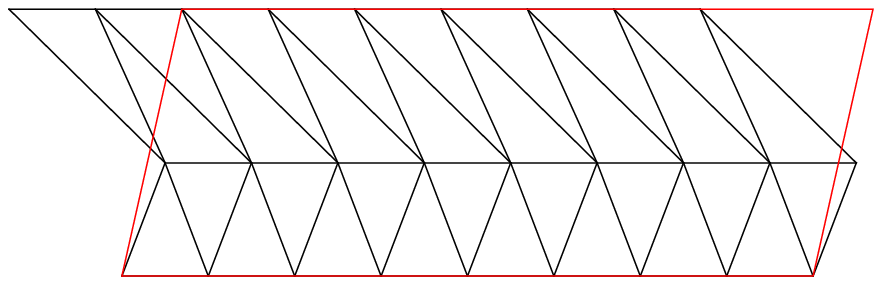}
\vskip1mm

$\DS-\frac{3}{16}-(-\frac{7}{16})=\frac{2}{8}$\hskip1.3cm
$\DS-\frac{2}{16}-(-\frac{6}{16})=\frac{2}{8}$\hskip1.5cm
$\DS-\frac{1}{16}-(-\frac{5}{16})=\frac{2}{8}$　

\vskip3mm

\includegraphics[height=1.3cm]{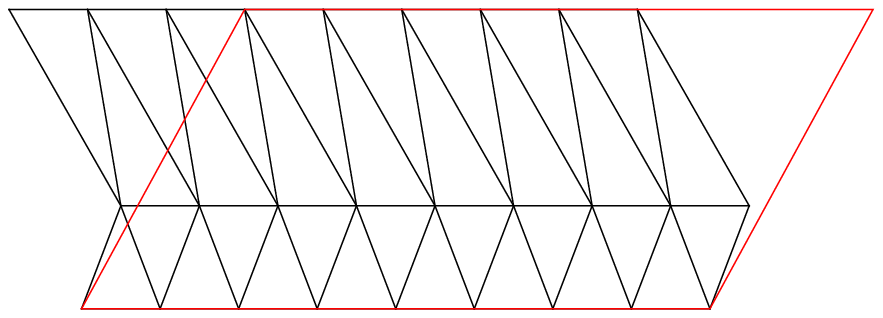}\ \ 
\includegraphics[height=1.3cm]{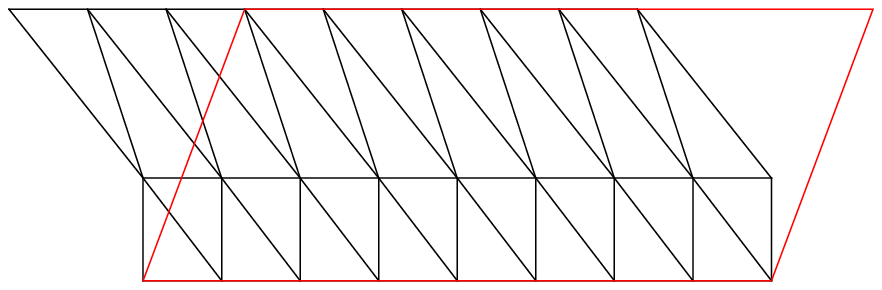}\ \ 
\includegraphics[height=1.3cm]{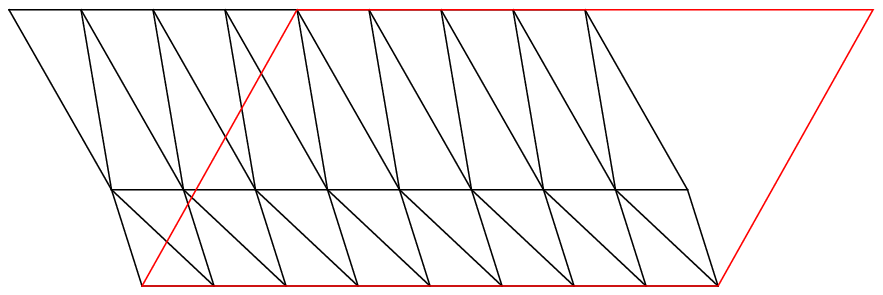}
\vskip1mm

$\DS-\frac{1}{16}-(-\frac{7}{16})=\frac{3}{8}$\hskip1.6cm
$\DS0-(-\frac{6}{16})=\frac{3}{8}$\hskip1.9cm
$\DS\frac{1}{18}-(-\frac{7}{16})=\frac{4}{8}$　

\end{center}
\caption{\small
Developments of flat tori $A_n^\rho\cup A_n^\sigma$ and $\sigma-\rho=\DS\frac{\ell}{n}$. 
Developments are seen from the side of the solid tori bounded by the flat tori. The upper parallelogram corresponds to the inner annulus $A_n^\rho$ and 
the lower parallelogram corresponds to the outer annulus $A_n^\sigma$. The common edges correspond to $P_0P_1\cdots P_7P_8$ and the upper edges and the lower edges
correspond to $Q_0^\sigma Q_1^\sigma \cdots Q_7^\sigma Q_8^\sigma$. Since $\sigma=\rho+\DS\frac{\ell}{n}$, the red parallelograms show the fundamental domains of the flat tori.
}\label{fig:rhosigma8dev}
\end{figure}

\section{Construction of embeddings of flat tori}\label{sect:tori}

In this section, we paste two suitable origami embedded annuli along boundaries and we construct 
origami embedded flat tori. 
Now we assume $n\geqq 5$. This construction is just a variant of that of Segerman (\cite{Segerman_video}, \cite{Segerman_book}, \cite{Segerman_shop}).

For two real numbers $\rho$, $\sigma\in (\DS\frac{1}{2},\DS\frac{1}{2}-\frac{1}{n})$, if origami embedded annuli $A_n^\rho$ and  $A_n^\sigma$ defined in Section \ref{sect:annuli} have the same boundaries and have disjoint inteiors, then we can paste $A_n^\rho$ and  $A_n^\sigma$ along the boundary and we obtain a surface which is an origami embedded flat torus. The fact that the obtained surface is a flat torus is easily verified because the development of the obtained surface is the union of two parallelograms which are the developments of  $A_n^\rho$ and  $A_n^\sigma$ pasted along the upper or lower edge. The identification of the edges of the development can be seen by the origami embedding and this determines the modulus of the origami embedded flat torus. 

If the difference of $\rho$ and $\sigma\in(-\DS\frac{1}{2},\DS\frac{1}{2}-\frac{1}{n})$ is a multiple of $\DS\frac{1}{n}$, the boundaries of $A_n^\rho$ coincides with those of $A_n^\sigma$. We would like to take $A_n^\rho$ inside of $A_n^\sigma$. 

For $\rho \in (-\DS\frac{1}{2},\DS-\frac{1}{2n}]$,  $A_n^\rho$ is placed closer to the real axis $\{0\}\times\RR$ than the one-sheet hyperboloid obtained by rotating the edge $P_0Q_1^\rho$ around the real axis $\{0\}\times\RR$
except the edges $P_kQ_{k+1}^\rho$ ($k=0$, \dots, $n-1$), and for $\rho \in [-\DS\frac{1}{2n},\DS\frac{1}{2}-\frac{1}{n})$, $A_n^\rho$ is placed closer to the real axis $\{0\}\times\RR$ than the one-sheet hyperboloid obtained by rotating the edge $P_0Q_0^\rho$ around the real axis $\{0\}\times\RR$ except the edges $P_kQ_k^\rho$ ($k=0$, \dots, $n-1$).

In order to insure the disjointness,
it is enough that the edges of $A_n^\sigma$ which are  closest to the real axis $\{0\}\times\RR$ is out side of this one-sheet hyperboloid.
Since the edge $P_0Q_0^\sigma$ is closer for $\sigma\in  (-\DS\frac{1}{2},-\frac{1}{2n}]$ and the edge 
$P_0Q_1^\sigma$ is closer for $\sigma\in  [-\DS\frac{1}{2n},\DS\frac{1}{2}-\frac{1}{n})$, we have the following conditions on $\rho$ and $\sigma$, which can be seen from the positions of the projections of edges on $\CC\times \{0\}$.
\begin{itemize}
\item
For $\rho \in (-\DS\frac{1}{2},\DS\frac{1}{2n}]$ and $
\sigma \in (-\DS\frac{1}{2},\DS\frac{1}{2n}]$,
$|\sigma|<|\rho+\DS\frac{1}{n}|$.
\item
For $\rho \in (-\DS\frac{1}{2},\DS\frac{1}{2n}]$ and $
\sigma \in [-\DS\frac{1}{2n},\DS\frac{1}{2}-\frac{1}{n})$,
$\sigma +\DS\frac{1}{n}<|\rho+\DS\frac{1}{n}|$.
\item
For $\rho \in
[-\DS\frac{1}{2n},\DS\frac{1}{2}-\frac{1}{n})$ and $
\sigma \in (-\DS\frac{1}{2},\DS\frac{1}{2n}]$,
$|\sigma|<|\rho|$.
\item 
For $\rho \in
[-\DS\frac{1}{2n},\DS\frac{1}{2}-\frac{1}{n})$ and $
\sigma \in [-\DS\frac{1}{2n},\DS\frac{1}{2}-\frac{1}{n})$,
$\sigma+\DS\frac{1}{n}<|\rho|$.
\end{itemize}
With the condition that 
$\sigma=\rho+\DS\frac{\ell}{n}$ ($\ell\in\ZZ$), the conditions are summed up as follows:
\begin{itemize}
\item
$\sigma=\rho+\DS\frac{\ell}{n}$ ($-\DS\frac{1}{2} <\rho<-\DS\frac{\ell}{2n}-\DS\frac{1}{n}$; $\ell=2$, \dots, $n-3$).
\item
$\sigma=\rho-\DS\frac{\ell}{n}$ ($\DS\frac{\ell}{2n}<\rho<\DS\frac{1}{2}-\DS\frac{1}{n}$; $\ell=2$, \dots, $n-3$).
\end{itemize}
See Figure \ref{fig:rhosigma}.

As concrete examples we construct 
origami embeddings $A_n^\rho\cup A_n^\sigma$  of flat tori by using a pair of the annuli shown in Figure \ref{fig:A_8_rho}.
 We have 12 pairs shown in Figure \ref{fig:rhosigma8} which satisfy the above  conditions.
The origami embeddings  
$A_n^\rho\cup A_n^\sigma$ are shown in Figure \ref{fig:rhosigma8emb}.
The developments of the origami embeddings $A_n^\rho\cup A_n^\sigma$
are shown in Figure \ref{fig:rhosigma8dev}.

\section{Moduli of embeddings of flat tori}\label{sect:moduli_of_tori}
In this section, by looking at the development of the surface $A_n^\rho\cup A_n^\sigma$, we compute the moduli of the flat tori and show the following theorem. 

\begin{theorem}\label{th:open_dense}
The construction of the origami embedding $A_n^\rho\cup A_n^\sigma$ for
$n\geqq 5$, $h>0$ and 
either  
$\sigma=\rho+\DS\frac{\ell}{n}$ ($-\DS\frac{1}{2} <\rho<-\DS\frac{\ell}{2n}-\DS\frac{1}{n}$; $\ell=2$, \dots, $n-3$) 
or
$\sigma =\rho-\DS\frac{\ell}{n}$ (
$\DS\frac{\ell}{2n}<\rho<\DS\frac{1}{2}-\DS\frac{1}{n}$; $\ell=2$, \dots, $n-3$) 
gives 
an origami embedding of a flat torus of any modulus except those represented by a pure imaginary number.
\end{theorem}

We describe the development of $A_n^\rho\cup A_n^\sigma$ more concretely.
We put in $\CC$ the vertex $\widehat{P}_0$ at the origin which corresponds to $P_0$. Set the $\CC$-coordinates of $\widehat{P}_k$  ($k=0$, \dots, $n$) by
$$\widehat{P}_k=2k\sin\DS\frac{\pi}{n}.$$ 
By using the computation on the point $X^\rho_1$ at the end of Section \ref{sect:annuli}, 
set the $\CC$-coordinates of $\ol{Q{}^\rho_k}$  ($k=0$, \dots, $n$)  by
$$\ol{Q{}^\rho_k}=
\sin(2\pi\rho+\DS\frac{\pi}{n})+(2k-1)\sin\DS\frac{\pi}{n}+\sqrt{-1}\sqrt{h^2+\big(\cos(2\pi\rho+\DS\frac{\pi}{n})-\cos\DS\frac{\pi}{n}\big)^2}.$$
Set the $\CC$-coordinates of $\ul{Q{}^\sigma_k}$  ($k=0$, \dots, $n$) by 
$$\ul{Q{}^\sigma_k}=
\sin(2\pi\sigma+\DS\frac{\pi}{n})+(2k-1)\sin\DS\frac{\pi}{n}-\sqrt{-1}\sqrt{h^2+\big(\cos(2\pi\sigma+\DS\frac{\pi}{n})-\cos\DS\frac{\pi}{n}\big)^2}.$$
Then the development of $A_n^\rho\cup A_n^\sigma$ is the union of
the parallelograms $\widehat{P}_0\widehat{P}_n \ol{Q_n^\rho}\,\ol{Q_0^\rho}$ 
and $\widehat{P}_0\widehat{P}_n \ul{Q_n^\sigma}\,\ul{Q_0^\sigma}$. 
The union of line segments $\ol{Q^\rho_0}\widehat{P}_0\ul{Q^\sigma_0}$ is
identified with $\ol{Q^\rho_n}\widehat{P}_n\ul{Q^\sigma_n}$.  
Since $\sigma=\rho\pm\DS\frac{\ell}{n}$, 
$\ul{Q^\sigma_{k}}\,\ul{Q^\sigma_{k+1}}$ is identified with
$\ol{Q^\rho_{k\pm\ell}}\,\ol{Q^\rho_{k+1\pm\ell}}$ ($k=0$, \dots, $n-1$ mod $n$).
Hence the modulus $\op{mod}(n,\pm \ell,\rho, h)=\DS\frac{\ol{Q_{\pm \ell}^\rho}-\ul{Q_0^\sigma}}{\widehat{P}_n}$ of $A_n^\rho\cup A_n^\sigma$ ($\sigma=\rho\pm\DS\frac{\ell}{n}$)
is computed as follows:
$$\ali
\op{mod}(n,\pm \ell,\rho, h)=&\DS\frac{\ol{Q_{\pm \ell}^\rho}-\ul{Q_0^\sigma}}{\widehat{P}_n}\\
=&
\DS\frac{1}{2n\sin\DS\frac{\pi}{n}}
\DS\bigg\{\sin(2\pi\rho+\DS\frac{\pi}{n})-\sin(2\pi\rho\pm\DS\frac{2\ell\pi}{n}+\DS\frac{\pi}{n})\pm 2\ell\sin\DS\frac{\pi}{n}
\\&\qquad\qquad\DS+\sqrt{-1}\bigg(\sqrt{h^2+\big(\cos(2\pi\rho+\DS\frac{\pi}{n})-\cos(\DS\frac{\pi}{n})\big)^2}\\&\qquad\qquad
+\sqrt{h^2+\big(\cos(2\pi\rho\pm\DS\frac{2\ell\pi}{n}+\DS\frac{\pi}{n})-\cos(\DS\frac{\pi}{n})\big)^2}\bigg)\bigg\}
\eali$$

We observe in this formula that the real part of the modulus does not depend on the height $h$. 
\begin{proposition}\label{prop:indep_of_height}
The real part of the modulus of the origami embedded flat torus 
$A_n^\rho\cup A_n^\sigma$ does not depend on the height $h$.
\end{proposition}

This proposition can also be shown  geometrically. For, when we deform the height $h$, the triangles $P_0P_1Q^\rho_1$ and $P_0Q_0^\rho Q^\rho_1$ change the height without changing the projected point $X_1^\rho$ on the line $P_0P_1$ of $Q^\rho_1$ and the projected point on the line  
 $Q_0^\rho Q^\rho_1$ of $P_0$, respectively. 

Taking account of Proposition \ref{prop:indep_of_height},
in order to determine the moduli of flat tori embeddable as $A_n^\rho\cup A_n^\sigma$, it is enough to look at the limit case where $h=0$. For, for a complex number of on the curve with $h=0$, the complex numbers with bigger imaginary part belong to the moduli.   

The curve with $h=0$ is given by 
$$\ali
\DS\frac{1}{2n\sin\DS\frac{\pi}{n}}
&\DS\bigg\{\sin(2\pi\rho+\DS\frac{\pi}{n})-\sin(2\pi\rho\pm\DS\frac{2\ell\pi}{n}+\DS\frac{\pi}{n})\pm 2\ell\sin\DS\frac{\pi}{n}\\
&\DS+\sqrt{-1}\bigg(\big|\cos(2\pi\rho+\DS\frac{\pi}{n})-\cos(\DS\frac{\pi}{n})\big|
+\big|\cos(2\pi\rho\pm\DS\frac{2\ell\pi}{n}+\DS\frac{\pi}{n})-\cos(\DS\frac{\pi}{n})\big|\bigg)\bigg\}.\eali$$
The cases where $\sigma=\rho+\DS\frac{\ell}{n}$ and
$\sigma=\rho-\DS\frac{\ell}{n}$ in Theorem \ref{th:open_dense} correspond to the cases 
where the real parts are positive and negative, respectively. They are the 
mirror images which have modulus symmetric with respect to the pure imaginary axis.
Hence to determine the moduli space for
$A_n^\rho\cup A_n^\sigma$, we look at only the cases where $\sigma-\rho=\DS\frac{\ell}{n}$ ($\ell=2$, \dots, $n-3$) and $\rho$ varies in $(-\DS\frac{1}{2},-\DS\frac{\ell}{2n}-\DS\frac{1}{n})$. 
Then the curve is
$$\ali
\DS\frac{1}{2n\sin\DS\frac{\pi}{n}}
&\DS\bigg\{\sin(2\pi\rho+\DS\frac{\pi}{n})-\sin(2\pi\rho+\DS\frac{2\ell\pi}{n}+\DS\frac{\pi}{n})+ 2\ell\sin\DS\frac{\pi}{n}\\
&\DS+\sqrt{-1}\bigg(\big|\cos(2\pi\rho+\DS\frac{\pi}{n})-\cos(\DS\frac{\pi}{n})\big|
+\big|\cos(2\pi\rho+\DS\frac{2\ell\pi}{n}+\DS\frac{\pi}{n})-\cos(\DS\frac{\pi}{n})\big|\bigg)\bigg\}.\eali$$
This curve is continuous on $\rho\in(-\DS\frac{1}{2},-\DS\frac{\ell}{2n}-\DS\frac{1}{n})$ for each $n$ and $\ell$ ($\ell=2$, \dots, $n-3$). 

We fix the ratio $\theta =\DS\frac{\ell}{n}$ and let $\ell$ and $n$ tend to the infinity and we see the curve converges to a rather comprehensible curve. That is, by putting $\ell_m=m\ell$, $n_m=mn$ and by letting tend $m$ to the infinity, the curve converges to 
$$c_\theta(\rho)=\DS\frac{1}{2\pi}\DS\big(\sin(2\pi\rho)-\sin(2\pi(\rho+\theta))+ 2\pi\theta\DS+\sqrt{-1}\bigg(2-\cos(2\pi\rho)-\cos(2\pi(\rho+\theta))\bigg)\big).$$
which is defined on $\rho\in (-\DS\frac{1}{2},-\DS\frac{\theta}{2})$ for rational $\theta\in(0,1)$.
\begin{proposition}
The closure of the union for all $\theta\in(0,1)$ of the images of curves $c_\theta(\rho)$ ($\rho\in (-\DS\frac{1}{2},-\DS\frac{\theta}{2})$) is the union of the following two domains: \\
the domain bounded by the line segment joining $0$ and $\DS\frac{\sqrt{-1}}{\pi}$ and the two cycloids 
$\{c_\theta(-\theta)\}_{\theta \in [0,\frac{1}{2}]}$ and $\{c_{\theta}(-\DS\frac{1}{2})\}_{\theta \in [0,\frac{1}{2}]}$;\\
the domain bounded by the three cycloids 
$\{c_{-2\rho}(\rho)\}_{\rho\in[-\frac{1}{2},0]}$,
$\{c_\theta(-\theta)\}_{\theta \in [0,\frac{1}{2}]}$ and
$\{c_{\theta}(-\DS\frac{1}{2})\}_{\theta \in [\frac{1}{2},1]}$.
\end{proposition}

\proof
We look at the complex valued function $\gamma(\rho,\theta)=c_\theta(\rho)$
defined on the interior of the triangle 
$\Delta$ with vertices $(\rho,\theta)=(-\DS\frac{1}{2},0)$, $(0,0)$, $(-\DS\frac{1}{2},1)$.
The partial derivatives are as follows:
$$\ali
\DS\frac{\pa \gamma}{\pa \rho}=&\cos(2\pi \rho)-\cos(2\pi(\rho+\theta))+\sqrt{-1}\big(\sin(2\pi\rho)+\sin(2\pi(\rho+\theta))\big)\\
\DS\frac{\pa \gamma}{\pa \theta}=&1-\cos(2\pi(\rho+\theta))+\sqrt{-1}\sin(2\pi(\rho+\theta))\eali
$$
The Jacobian of $\gamma:\Delta\lra \CC$ is calculated as follows:
$$\ali
&\{\cos(2\pi \rho)-\cos(2\pi(\rho+\theta))\}\sin(2\pi(\rho+\theta)
\\&-\{\sin(2\pi\rho)+\sin(2\pi(\rho+\theta))\}\{1-\cos(2\pi(\rho+\theta))\}
\\=&
\cos(2\pi \rho)\sin(2\pi(\rho+\theta))+
\sin(2\pi \rho)\cos(2\pi(\rho+\theta))
-\{\sin(2\pi\rho)+\sin(2\pi(\rho+\theta))\}\\=&
\sin(2\pi(2\rho+\theta))-
\sin(2\pi\rho)-\sin(2\pi(\rho+\theta))\\
=&4\sin(\pi \rho)\sin(\pi(\rho+\theta))\sin(\pi(2\rho+\theta))
\eali$$
Thus in the interior of the triangle $\Delta$, the map is singular along the line $\rho+\theta=0$. In fact, along the line $\rho+\theta=0$, $\DS\frac{\pa \gamma}{\pa \rho}=\cos(2\pi\rho)-1+\sqrt{-1}\,\sin(2\pi\rho)$ (the real part and the imaginary part are both negative for $\rho\in (\DS-\frac{1}{2},0)$) and $\DS\frac{\pa \gamma}{\pa \theta}=0$.
This line corresponds to the place where the limit of $\sigma=\rho+\theta$ is zero and the image of this line is given by the following curve 
$c_\theta(-\theta)$:
$$c_\theta(-\theta)=\DS\frac{2\pi\theta-\sin(2\pi\theta)+ \sqrt{-1}(1-\cos(2\pi\theta)))}{2\pi}\quad (\theta \in (0,\DS\frac{1}{2})),
$$
which is a cycloid joining $0$ and $\DS\frac{1}{2}+\DS\frac{\sqrt{-1}}{2\pi}\in\CC$.

\begin{figure}
\begin{center}
\includegraphics[height=5cm]{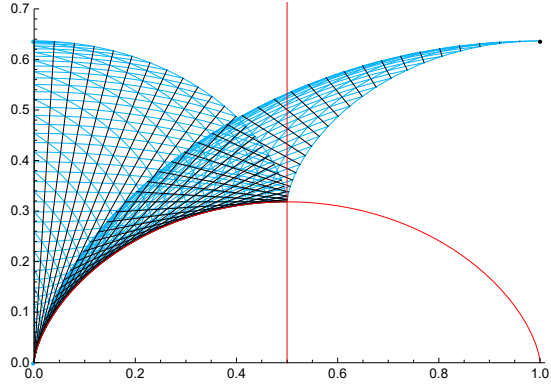}
\end{center}
\caption{\small
The image in $\CC$ of $c_\theta(\rho)$ ($(\rho,\theta)\in\Delta$). For fixed $\theta$, the curve $c_\theta(\rho)$ is shown in black.}\label{fig:Image_c_theta_rho}
\end{figure}

The images of the edges of $\Delta$ is as follows:
\begin{itemize}
\item
The edge joining $(0,0)$ and $\DS(-\frac{1}{2},1)$ where $\theta =-2\rho$ is mapped to the following curve $c_{-2\rho}(\rho)$: 
$$c_{-2\rho}(\rho)=\DS\frac{2\sin(2\pi\rho)-4\pi\rho+\sqrt{-1}(2-2\cos(2\pi\rho))}{2\pi}\quad (\rho\in (-\DS\frac{1}{2},0))$$
which is again a cycloid joins $0$ and $1+\DS\frac{2}{\pi}\sqrt{-1}$.
\item
The edge joining  $\DS(-\frac{1}{2},1)$ and $\DS(-\frac{1}{2},0)$ where $\rho=-\DS\frac{1}{2}$ is mapped to the following curve $c_{\theta}(-\DS\frac{1}{2})$:
$$c_{\theta}(-\DS\frac{1}{2})=\DS\frac{\sin(2\pi\theta)+2\pi\theta+\sqrt{-1}(3+\cos(2\pi\theta))}{2\pi}\quad(\theta\in(0,1)).$$
This is a cycloid, and 
$$\op{Im}(c_{\theta}(-\DS\frac{1}{2}))=\op{Im}(c_{1-\theta}(-\DS\frac{1}{2}))
\ \ \text{and}\ \ \op{Re}(c_{\theta}(-\DS\frac{1}{2}))=1-\op{Re}(c_{1-\theta}(-\DS\frac{1}{2})),$$ 
i.e., the image is symmetric with respect to the line $\op{Re}(z)=\DS\frac{1}{2}$. The point $\theta=\DS\frac{1}{2}$ is a critical point whose image is a cusp. The curve $c_{\theta}(-\DS\frac{1}{2})$ ($\theta\in(0,1)$) joins $\DS\frac{2}{\pi}\sqrt{-1}$, $\DS\frac{1}{2}+ \DS\frac{\sqrt{-1}}{\pi}$ and $1+\DS\frac{2}{\pi}\sqrt{-1}$.
\item
The edge joining  $\DS(-\frac{1}{2},0)$ and $(0,0)$ where $\theta=0$ is mapped to the following curve $c_0(\rho)$: 
$$c_0(\rho)=\sqrt{-1}\DS\frac{1-\cos(2\pi\rho)}{\pi}\quad(\rho\in(-\DS\frac{1}{2},0))$$
which is a line segment joining $0$ and $\DS\frac{2}{\pi}\sqrt{-1}$. See Figure \ref{fig:Image_c_theta_rho}. 
\end{itemize}

The above curves gives the boundary simple closed curves which are boundaries of the domains in the proposition. Since the map is regular in the interior of $\Delta\setminus \{\rho+\theta=0\}$, the image of each component of 
$\Delta\setminus \{\rho+\theta=0\}$ coincides with the interior of one of the two domains in the proposition.
Thus the proposition is shown. \hfill $\square$
\vskip5mm

\begin{remark}\label{remark:line_segment}
Since $$\ali
\DS\frac{\pa \gamma}{\pa \rho}=&\cos(2\pi \rho)-\cos(2\pi(\rho+\theta))+\sqrt{-1}\big(\sin(2\pi\rho)+\sin(2\pi(\rho+\theta))\big)\\=&
2\sin\pi(2\rho+\theta)\sin(\pi\theta)
+2\sqrt{-1}\sin\pi(2\rho+\theta)\cos(\pi\theta)\eali$$
and the ratio $\op{Im}(\DS\frac{\pa\gamma}{\pa \rho})\big/\op{Re}(\DS\frac{\pa\gamma}{\pa \rho})$
of the imaginary part and the real part is  
$\DS\cot(\pi\theta)$ which is independent of $\rho$, $\{c_\theta(\rho)\}_{\rho\in (-\frac{1}{2},-\frac{\theta}{2})}$ is a line segment joining the two points $c_\theta(-\DS\frac{1}{2})$ of the cycloid $\{c_\theta(-\DS\frac{1}{2})\}_{\theta \in (\frac{1}{2},1)}$ and $c_\theta(-\DS\frac{\theta}{2})$ of the 
cycloid $\{c_{-2\rho}(\rho)\}_{\rho \in (-\frac{1}{2}, 0)}$.
For $0<\theta <\DS\frac{1}{2}$, the line segment $\{c_\theta(\rho)\}_{\rho\in (-\frac{1}{2},-\frac{\theta}{2})}$ is tangent to the cycloid $\{c_\theta(-\theta)\}_{\theta \in (0,\frac{1}{2})}$, and
for $\DS\frac{1}{2}<\theta <1$, the prolongation of the line segment is tangent to the same cycloid.\end{remark}

\noindent 
\textit{Proof} of Theorem \ref{th:open_dense}.
We take any rational number $\theta$  in the open interval $(0,\DS\frac{1}{2})$.
If $\rho$ is close to $-\theta\in (-\DS\frac{1}{2},-\DS\frac{\theta}{2})$, then the point $c_\theta(\rho)$ is on the tangent line of the cycloid $\{c_\theta(-\theta)\}_{\theta \in (0,\frac{1}{2})}$ at the point $c_\theta(-\theta)$ and close to the point $c_\theta(-\theta)$ (See Remark \ref{remark:line_segment}). Then the modulus $\op{mod}(n,\ell,\rho,h)$ for $\theta =\DS\frac{\ell}{n}$ and small positive $h$ stays near the tangent line $c_\theta(\rho)$. Then the points with greater imaginary part than $\op{mod}(n,\ell,\rho,h)$ are moduli
of flat tori which can be origami embedded in $\RR^3$, and so do the points with greater imaginary part than the cycloid $\theta\longmapsto c_\theta(-\theta)$. 

Note that the real part $\op{Re}(
c_\theta(-\theta))=\DS\frac{2\pi\theta-\sin(2\pi\theta)}{2\pi}$ and
the imaginary part $\op{Im}(c_\theta(-\theta))=\DS\frac{1-\cos(2\pi\theta))}{2\pi}$
of $c_\theta(-\theta)$ are both monotone increasing with respect to $\theta\in (0,\DS\frac{1}{2})$ 
and $\DS\lim_{\theta\to0} c_\theta(-\theta)=0\in\CC$ and $\DS\lim_{\theta\to1/2} c_\theta(-\theta)=\DS\frac{1}{2}+\DS\frac{i}{2\pi}\in\CC$.
Since $\op{Im}(c(\theta))<\DS\frac{1}{2\pi}<\DS\frac{\sqrt{3}}{2}$, the curve $\{c(\theta)\ \big| \theta\in(0,\DS\frac{1}{2})\}$ is contained in the unit disk in $\CC$. This implies that the set of moduli of the origami embedded flat tori $A_n^\rho\cup A_n^\sigma$ ($n\geqq5$; $h>0$, $\sigma=\rho+\DS\frac{\ell}{n}$;  $\ell=2$, \dots, $n-3$) contains the set $$\{x+y\sqrt{-1}\ \big|\ x\in (0,\DS\frac{1}{2}),\ y\geqq \sqrt{1-x^2}\}$$
which is almost the half of the fundamental domain of moduli of flat tori.

By looking at the case where $\sigma=\rho-\DS\frac{\ell}{n}$, we obtain the moduli which are mirror images of those in the case where 
 $\sigma=\rho+\DS\frac{\ell}{n}$, and hence we showed Theorem \ref{th:open_dense} except for the moduli with the real part $\DS\frac{1}{2}$.

To treat the case of the real part $\DS\frac{1}{2}$, 
we look at the line segment  $\{c_\theta(\rho)\}_{\rho\in(-\frac{1}{2},-\frac{\theta}{2})}$ for $1>\theta>\DS\frac{1}{2}$ (See Remark \ref{remark:line_segment}).
As we remarked, $c_\theta(\rho)$ connects the point $c_\theta(-\DS\frac{1}{2})$ of the cycloid $\{c_\theta(-\DS\frac{1}{2})\}_{\theta \in (\frac{1}{2},1)}$ and the point  $c_\theta(-\DS\frac{\theta}{2})$ of the 
cycloid $\{c_{-2\rho}(\rho)\}_{\rho \in (-\frac{1}{2}, -\frac{1}{4})}$.
Thus the curve $c_\theta(\rho)$ for $\theta >\DS\frac{1}{2}$ close to $\DS\frac{1}{2}$ crosses the line where the real part is $\DS\frac{1}{2}$ near the point $\DS\frac{1}{2}+\frac{\sqrt{-1}}{\pi}$. Hence by the argument as before, the line $\DS\frac{1}{2}+y\sqrt{-1}$ ($y> \DS\frac{\sqrt{3}}{2}$) is moduli of origami embedded flat tori. 
\hfill $\square$
\vskip5mm

\begin{remark}
In order to represent moduli near the imaginary axis, we need to use $A_n^\rho\cup A_n^\sigma$ ($\sigma=\rho\pm\DS\frac{\ell}{n}$) with large $n$.
\end{remark}

\section{Origami embeddings of flat tori of pure imaginary moduli}\label{sect:rectangular}

In order to treat the flat tori with pure imaginary moduli, we use the following simple observation.
\begin{proposition}\label{prop:cut_along_z=a}
For the emmbedded annulus $A^\rho_n\subset \CC\times \RR$ of height $h$ defined in Section \ref{sect:annuli} and $0<a<h$, consider the parts $A^\rho_n|_{z\in[0,a]}$ and $A^\rho_n|_{z\in[a,h]}$ of $A^\rho_n$ where the $\RR$ coordinate $z$ belong to $[0,a]$ and to $[a,h]$, respectively.
Then the developments of  $A^\rho_n|_{z\in[0,a]}$ and $A^\rho_n|_{z\in[a,h]}$ are parallelograms 
obtained as the parts lower and upper than the horizontal line dividing the development of $A^\rho_n$ in the ratio $a:h-a$.
\end{proposition}

\proof 
It follows from the fact that the intersection of triangles $\triangle P_0P_1Q_1^\rho$, 
$\triangle Q_1^\rho Q_0^\rho P_0$ and the plane $\{z=a\}$ are the line segments pararell to 
$P_0P_1$, $Q_1^\rho Q_0^\rho$. \hfill $\square$
\vskip5mm

The following proposition completes the proof of Theorem \ref{th:main}.

\begin{proposition}\label{prop:rectangular}
The flat tori with rectangular fundamental domains can be origami embedded in the 3-dimensional Euclidean space.
\end{proposition}

\proof
We take an origami embedding $A_n^\rho\cup A_n^\sigma$  
of a flat torus constructed in Section \ref{sect:tori}. 
We cut $A_n^\rho\cup A_n^\sigma$ along the plane $\{z=a\}$ for $0<a<h$ and
we take the double of $A_n^\rho|_{z\in[0,a]}\cup A_n^\sigma|_{z\in[0,a]}$ or
$A_n^\rho|_{z\in[a,h]}\cup A_n^\sigma|_{z\in[a,h]}$.
Then the double has a rectangular fundamental domain. The modulus is small if the height of the double is small, and hence any point of the imaginary axis can be realized as an origami embedded 
flat torus.
\hfill $\square$
\vskip5mm

\begin{figure}
\begin{center}
\includegraphics[height=3cm]{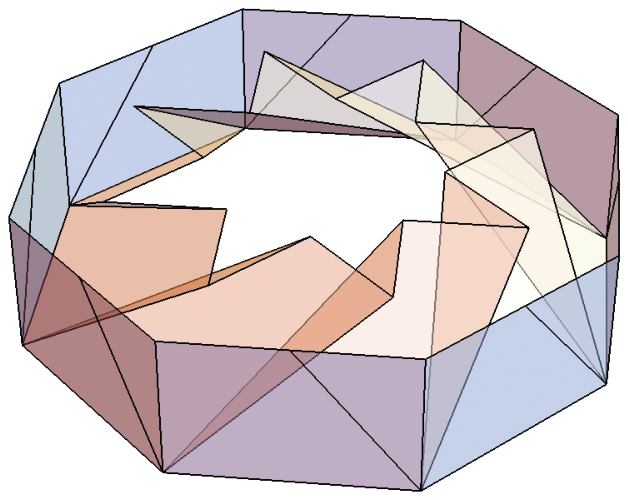}
\includegraphics[height=3.6cm]{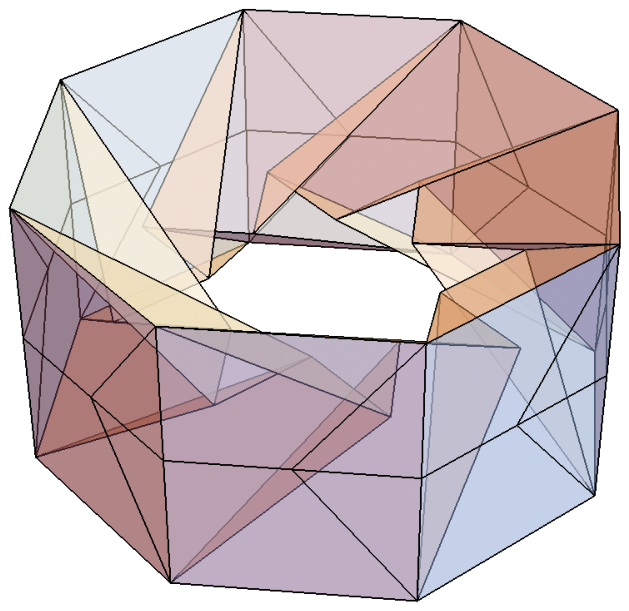}

\hskip20mm\includegraphics[height=1.8cm]{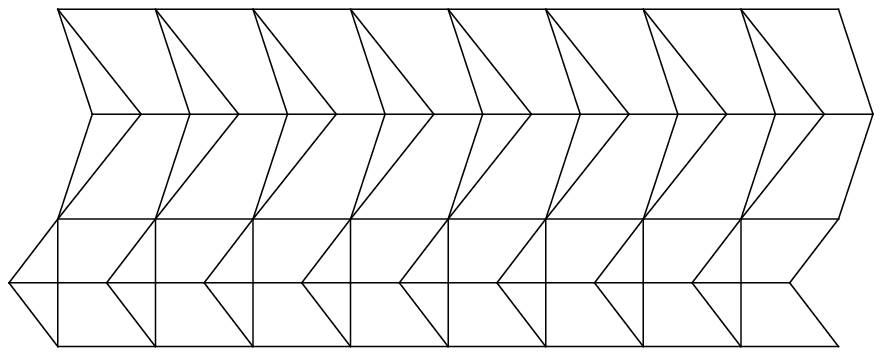}
\vskip3mm

\includegraphics[height=3cm]{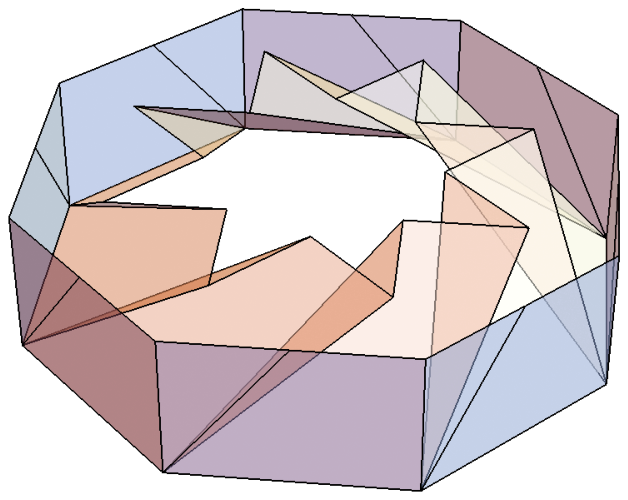}
\includegraphics[height=3.6cm]{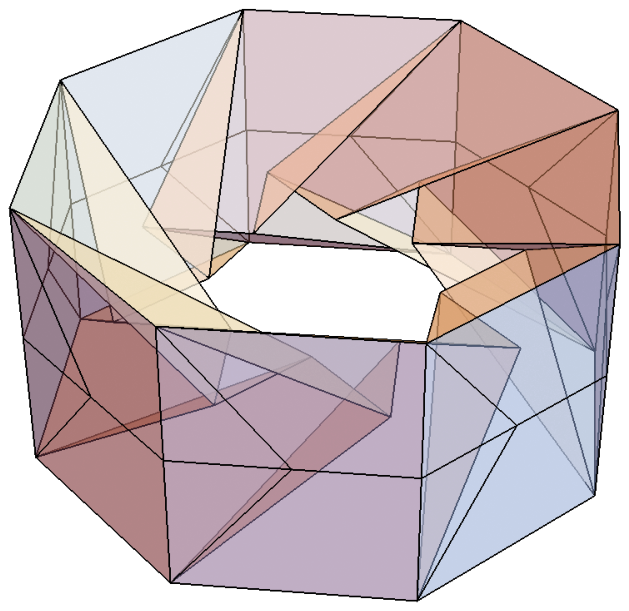}

\hskip20mm\includegraphics[height=1.8cm]{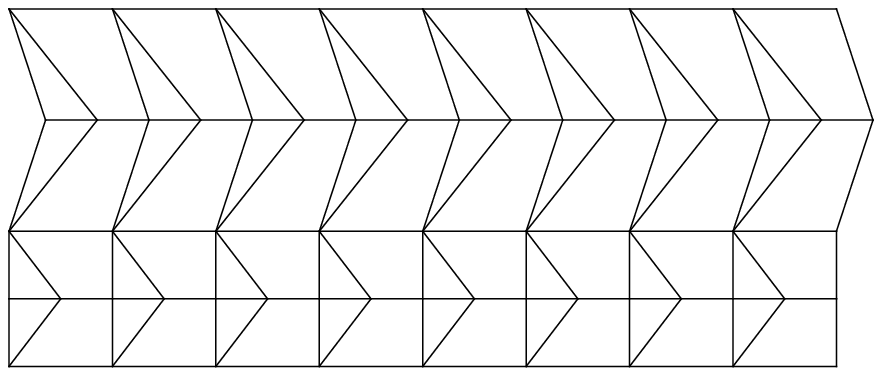}

\caption{\small $A_n^\rho|_{z\in[0,\frac{h}{2}]}\cup A_n^\sigma|_{z\in[0,\frac{h}{2}]}$ 
for $n=8$, $h=1$, $(\rho,\sigma)=(\DS\frac{4\pi}{8},-\DS\frac{2\pi}{8})$, $(\DS\frac{4\pi}{8},0)$, 
the doubles of them and the development of the doubles.}\label{fig:half_double_dev}
\end{center}
\end{figure}

We give explicit examples of origami embedded flat tori with rectangular fundamental domain in Figure \ref{fig:half_double_dev}.

\begin{remark}
Proposition \ref{prop:cut_along_z=a} can be generalized so that for $0<a<b<h$, the development of $A_n^\rho|_{z\in[a,b]}$ is a parallelogram and one can try to use it to construct some other origami embeddings of flat tori. However we think it is necessary to contain a regular $n$-gon in our constructions until now. The possible reason is as follows: The intersection $A_n^\rho\cap \{z=a\}$ is a $2n$-gon and from the shape of this intersection we find $n$, $\rho$ and $\DS\frac{a}{h}$ except the case where the intersection is regular $2n$-gon. That is, if the intersection is not a regular $2n$-gon, $n$ is obtained as the half of the number of vertices, the sum of the lengths of edges is eaual to that of the regular $n$-gon in the construction, the ratio of the lengths of adjacent edges gives $\DS\frac{a}{h-a}$, and the angle of the edges of the intersection gives $\rho$ because they are parallel to edges of the top or bottom regular $n$-gon. 
Thus if there are no cross sections which are regular $n$-gons,  
we can only use $A_n^\rho|_{z\in[a,b]}$ with fixed $\rho$, and only give rise to annuli.
\end{remark}

\begin{remark}
The doubling construction of this section is essentially the same as the 
bending construction for triangular cylinders given by Zalgaller \cite{Zalgaller}. In view of Proposition \ref{prop:indep_of_height}, our Proposition \ref{prop:rectangular} follows from the construction by Zalgaller \cite{Zalgaller}.
\end{remark}

\end{document}